\documentclass[a4paper,12pt,twoside]{article}
\usepackage{amsmath,amssymb,ifthen}
\usepackage[amsmath,thmmarks]{ntheorem}
\usepackage{paper}
\usepackage[all]{xy}
\usepackage{mathrsfs}
\newcommand{\red}{\mathrm{red}}
\newcommand{\ad}{\mathrm{ad}}
\DeclareMathOperator{\Spf}{Spf}
\DeclareMathOperator{\Spa}{Spa}
\DeclareMathOperator{\cosk}{cosk}
\DeclareMathOperator{\Adj}{Adj}
\SelectTips{cm}{11}
\title{Variants of formal nearby cycles}
\author{Yoichi Mieda}
\begin{document}
\maketitle

\begin{firstfootnote}
 Faculty of Mathematics, Kyushu University, 744 Motooka, Nishi-ku, Fukuoka, 819--0395 Japan

 E-mail address: \texttt{mieda@math.kyushu-u.ac.jp}

\end{firstfootnote}

\begin{abstract}
 In this paper, we introduce variants of formal nearby cycles for a locally noetherian formal scheme
 over a complete discrete valuation ring. If the formal scheme is locally algebraizable, 
 then our nearby cycle gives a generalization of Berkovich's formal nearby cycle.
 Our construction is entirely scheme-theoretic and does not require rigid geometry.
 Our theory is intended for applications to the local study of the cohomology of Rapoport-Zink spaces.
\end{abstract}

\section{Introduction}
In the papers \cite{MR1262943} and \cite{MR1395723}, Berkovich defined the formal nearby cycle functor for 
a formal scheme $\mathcal{X}$ over a complete discrete valuation ring $R$, and proved the comparison result
that if $\mathcal{X}$ is obtained by completing a scheme $X$ locally of finite type over $R$
along a closed subscheme $Y$ of the special fiber $X_s$ of $X$, then the formal nearby cycle is isomorphic to
the restriction of the nearby cycle $R\psi_X$ of $X$ to $Y$.
In particular, $(R\psi_X\Lambda)\vert_Y$ depends only on the completion of $X$ along $Y$ 
for $\Lambda=\Z/\ell^n\Z$.
The theory of Berkovich plays an important role in the study of the cohomology of 
Shimura varieties and Rapoport-Zink spaces; for example, see \cite{MR1876802} and \cite{MR2074714}.

In this paper, we introduce variants of Berkovich's formal nearby cycle and investigate their properties.
For simplicity, assume that the residue field of $R$ is separably closed.
Our results are roughly summarized in the following theorem:

\begin{thm}\label{thm:main-thm}
 For a locally noetherian excellent formal scheme $\mathcal{X}$ over $\Spf R$ such that $\mathcal{X}_\red$ is
 separated, and a pair 
 $\mathcal{Z}=(\mathcal{Z}_1,\mathcal{Z}_2)$ of closed formal subschemes of the special fiber $\mathcal{X}_s$
 of $\mathcal{X}$,
 we can define the object $R\Psi_{\!\mathcal{X},\mathcal{Z}}\Lambda$ of 
 $D^+(\mathcal{X}_\red,\Lambda)$ (for the definition of excellent schemes, see Section \ref{subsec:formal}).
 It enjoys the following properties:
 \begin{enumerate}
  \item $R\Psi_{\!\mathcal{X},\mathcal{Z}}\Lambda$ is functorial on $\mathcal{X}$.
	In particular, if a group $G$ acts on $\mathcal{X}$ and $\mathcal{Z}$ is stable under the action,
	then $R\Psi_{\!\mathcal{X},\mathcal{Z}}\Lambda$ has a natural $G$-equivariant structure.
  \item If $\mathcal{X}$ is obtained by completing a scheme $X$ locally of finite type over $R$
	along a closed subscheme $Y$ of $X_s$ and $\mathcal{Z}$ comes from
	a pair $Z=(Z_1,Z_2)$ of closed subschemes of $X_s$,
	then we have a functorial isomorphism
	\[
	R\Psi_{\!\mathcal{X},\mathcal{Z}}\Lambda\cong (Rj_*Rj^!R\psi_X\Lambda)\vert_Y,
	\]
	where $j$ denotes the natural immersion $Z_2\setminus Z_1\hooklongrightarrow X_s$.
  \item If $\mathcal{X}$ is locally algebraizable (Definition \ref{defn:loc-alg}) or adic over $\Spf R$,
	then $R\Psi_{\!\mathcal{X}}\Lambda:=R\Psi_{\!\mathcal{X},(\varnothing,\mathcal{X}_s)}\Lambda$ is
	canonically isomorphic to Berkovich's formal nearby cycle for $\mathcal{X}$.
  \item If $\mathcal{X}$ is quasi-compact, separated, pseudo-compactifiable (Definition \ref{defn:pseudo-comp})
	and locally algebraizable, we have the following
	natural isomorphism 
	\[
	  H^q_c\bigl(t(\mathcal{X})_{\overline{\eta}},\Lambda\bigr)\cong H^q_c(\mathcal{X}_\red,R\Psi_{\!\mathcal{X},c}\Lambda),
	\]
	where $t(\mathcal{X})_{\overline{\eta}}$ denotes the geometric generic fiber of $\mathcal{X}$ and 
	$R\Psi_{\!\mathcal{X},c}\Lambda$ denotes $R\Psi_{\!\mathcal{X},(\varnothing,\mathcal{X}_\red)}\Lambda$.
 \end{enumerate}
\end{thm}
The property ii) in the theorem above says that our nearby cycle preserves the information of $R\psi_X\Lambda$
outside $Y$, which is discarded if we consider Berkovich's formal nearby cycle.
The information outside $Y$ seems important in the study of the cohomology of Rapoport-Zink spaces.
For example, the content of this article is used in an essential manner in the paper \cite{RZ-GSp4}, 
where we consider the non-cuspidality of the cohomology of the Rapoport-Zink space for $\mathrm{GSp}(4)$.
See also \cite{non-cusp}.

The property iv) is also remarkable. 
This is a new phenomenon; note that the compactly supported cohomology of Berkovich's formal nearby cycle
is not necessarily equal to the compactly supported cohomology of $t(\mathcal{X})_{\overline{\eta}}$
(\cf Remark \ref{rem:xi-isom}).

Our method of constructing $R\Psi_{\!\mathcal{X},\mathcal{Z}}\Lambda$ is completely different
from Berkovich's one. 
Berkovich used rigid geometry (or more precisely, his own theory of analytic spaces),
while we only use the scheme theory.
For an affine formal scheme $\Spf A$, 
we consider the nearby cycle of the affine scheme $\Spec A$ (with some modification).
For a general formal scheme, we patch the local construction above by using the simplicial technique.

We sketch the outline of the paper. In Section 2, we collect some facts which are required to construct
$R\Psi_{\!\mathcal{X},\mathcal{Z}}\Lambda$. In particular, the content in Section \ref{subsec:formal}
is crucial for our construction.
In Section 3, we give a definition of our nearby cycle
$R\Psi_{\!\mathcal{X},\mathcal{Z}}\Lambda$ and prove various properties on it.
In Section 4, we compare our theory with the theory of rigid spaces.
We use the framework of adic spaces due to Huber. First we recall the definition of Berkovich's formal
nearby cycle functor and prove the comparison result Theorem \ref{thm:main-thm} iii).
Next we will give Theorem \ref{thm:main-thm} iv), whose proof is rather involving.

\bigbreak

\noindent{\bfseries Acknowledgment}\quad
The author would like to thank Tomoyuki Abe and Atsushi Shiho for the stimulating discussions.

\bigbreak

\noindent{\bfseries Notation}\quad
Let $R$ be a complete discrete valuation ring with the separably closed residue field $k$,
$F$ the fraction field of $R$ and $\varpi$ a uniformizer of $R$. Put $S=\Spec R$, $\mathcal{S}=\Spf R$. 
We fix a separable closure $\overline{F}$ of $F$ and denote the closed 
(resp.\ generic, resp.\ geometric generic) point of $S$ by $s$
(resp.\ $\eta$, resp.\ $\overline{\eta}$). We also regard $s$ as a closed subscheme of $\mathcal{S}$.
We fix a prime number $\ell$ which is invertible in $R$, and set $\Lambda=\Z/\ell^n\Z$ with $n\ge 0$.
For an $S$-scheme $X$, $R\psi_X$ denotes the nearby cycle functor for $X\longrightarrow S$.

Every sheaf and cohomology are considered in the \'etale topology.

\section{Preliminaries}\label{sec:preliminaries}
In this section, we give some preparations for constructing our nearby cycle 
$R\Psi_{\!\mathcal{X},\mathcal{Z}}\Lambda$.

\subsection{Preliminaries on open coverings}\label{subsec:open-cov}
Let $X$ be a scheme. For an open covering $U=\{U_i\}_{i\in I}$ of $X$, we put $U_0=\coprod_{i\in I}U_i$
and let $a\colon U_\bullet\longrightarrow X$ be the augmented simplicial scheme $\cosk_0(U_0/X)$.
Concretely, $U_m$ is the $(m+1)$-fold fiber product of $U_0$ over $X$ for every $m\ge 0$.
In particular, $U_m$ is a disjoint union of open subschemes of $X$.
We call $a\colon U_\bullet\longrightarrow X$ the hypercovering associated with $U$.

\begin{defn}
 Let $\mathcal{F}=(\mathcal{F}^m)_{m\ge 0}$ be a $\Lambda$-sheaf on $U_\bullet$. We say that
 $\mathcal{F}$ is cartesian if for every structure morphism $\phi\colon U_m\longrightarrow U_n$ of $U_\bullet$,
 $\phi^*\mathcal{F}^n\longrightarrow \mathcal{F}^m$ is an isomorphism.
 We denote by $D^+_{\mathrm{cart}}(U_\bullet,\Lambda)$ the full subcategory of
 $D^+(U_\bullet,\Lambda)$ consisting of lower bounded complexes whose cohomology are all cartesian.
\end{defn}

The following lemma is clear:

\begin{lem}\label{lem:cartesian-sheaf}
 A $\Lambda$-sheaf $\mathcal{F}=(\mathcal{F}^m)_{m\ge 0}$ on $U_\bullet$ is cartesian
 if and only if there exists a sheaf $\mathcal{G}$ on $X$ such that $\mathcal{F}\cong a^*\mathcal{G}$.
\end{lem}

\begin{prop}\label{prop:cart-complex}
 The functor $a^*\colon D^+(X,\Lambda)\longrightarrow D^+_{\mathrm{cart}}(U_\bullet,\Lambda)$ gives a
 categorical equivalence. The functor $Ra_*$ gives the quasi-inverse of $a^*$.
\end{prop}

\begin{prf}
 First note that $Ra_*\circ a^*\cong \id$, since $a$ is a morphism of cohomological descent.
 
 For an object $L$ of $D^+_{\mathrm{cart}}(U_\bullet,\Lambda)$,
 we will prove that $a^*Ra_*L\longrightarrow L$ is an isomorphism.
 This is equivalent to that $a^*R^ma_*L\longrightarrow H^m(L)$ is an isomorphism for every integer $m$.
 We fix $m$ and put $L'=\tau_{\le m}L$. Then, by the distinguished triangle
 $L'\longrightarrow L\longrightarrow \tau_{>m}L\longrightarrow L'[1]$,
 the morphism $a^*R^ma_*L'\longrightarrow a^*R^ma_*L$ is an isomorphism.
 On the other hand, $H^m(L')\longrightarrow H^m(L)$ is obviously an isomorphism.
 Therefore we may replace $L$ by $L'$. 
 In other words, we are reduced to the case where $L$ is bounded.
 Moreover, we may assume that $H^i(L)=0$ for $i<0$.

 Under these conditions, we prove the isomorphy of $a^*Ra_*L\longrightarrow L$ by 
 the induction on $d=\max\{i\in\Z\mid H^i(L)\neq 0\}$.
 If $d=0$, by Lemma \ref{lem:cartesian-sheaf}, there exists a $\Lambda$-sheaf $\mathcal{F}$ on $X$
 such that $L\cong a^*\mathcal{F}$. We can replace $L$ by $a^*\mathcal{F}$.
 Since the composite of adjunction morphisms 
 $a^*\mathcal{F}\yrightarrow{\cong} a^*Ra_*a^*\mathcal{F}\longrightarrow a^*\mathcal{F}$ is the identity,
 $a^*Ra_*a^*\mathcal{F}\longrightarrow a^*\mathcal{F}$ is an isomorphism. 
 Note that the first morphism is an isomorphism since $a$ is a morphism of cohomological descent.
 For a general $d$, consider the following commutative diagram whose rows are distinguished triangles: 
 \[
  \xymatrix{%
 a^*Ra_*(\tau_{\le d-1}L)\ar[r]\ar[d]^-{(1)}& a^*Ra_*L\ar[r]\ar[d]^-{(2)}& a^*Ra_* H^d(L)[-d]\ar[r]\ar[d]^-{(3)}& a^*Ra_*(\tau_{\le d-1}L)[1]\ar[d]^-{(1)}\\
 \tau_{\le d-1}L\ar[r]& L\ar[r]& H^d(L)[-d]\ar[r]& \tau_{\le d-1}L[1]\lefteqn{.}
 }
 \]
 The morphism (1) is an isomorphism by the induction hypothesis. 
 We have already seen that the morphism (3) is an isomorphism.
 Thus the morphism (2) is also an isomorphism.
 This completes the proof.
\end{prf}

Let $V=\{V_j\}_{j\in J}$ be another open covering of $X$ and $b\colon V_\bullet\longrightarrow X$ the associated
hypercovering. From $a\colon U_\bullet\longrightarrow X$ and $b\colon V_\bullet\longrightarrow X$, we can construct
the augmented bisimplicial scheme $c\colon W_{\bullet\bullet}\longrightarrow X$, where
$W_{mn}=U_m\times_XV_n$. 

\begin{lem}\label{lem:bisimplicial-descent}
 The functor $c^*\colon D^+(X,\Lambda)\longrightarrow D^+(W_{\bullet\bullet},\Lambda)$ is fully faithful.
\end{lem}

\begin{prf}
 We have the canonical augmentation $p\colon W_{\bullet\bullet}\longrightarrow U_{\bullet}$
 such that $c=a\circ p$ and the functor $c^*$ is the composite of
 $D^+(X,\Lambda)\yrightarrow{a^*} D^+(U_\bullet,\Lambda)\yrightarrow{p^*} D^+(W_{\bullet\bullet},\Lambda)$.
 Since we already know that $a^*$ is fully faithful, it suffices to show that $p^*$ is fully faithful,
 or equivalently, $Rp_*p^*\cong \id$.

 Let $L$ be an object of $D^+(U_\bullet,\Lambda)$. To prove that $L\longrightarrow Rp_*p^*L$ is an isomorphism,
 it is sufficient to prove that its restriction to $U_m$ is an isomorphism for every $m\ge 0$.
 The restriction obviously coincides with the adjunction morphism
 $L_m\longrightarrow Rp_{m*}p_m^*L_m$, where $L_m=L\vert_{U_m}$ and $p_m\colon W_{m\bullet}\longrightarrow U_m$
 is the augmentation induced from $p$.
 Since $p_m$ is the base change of $b\colon V_\bullet\longrightarrow X$ by $U_m\longrightarrow X$, 
 it is a morphism of cohomological descent. Thus $L_m\longrightarrow Rp_{m*}p_m^*L_m$ is an isomorphism,
 as desired.
\end{prf}

\subsection{Preliminaries on formal schemes}\label{subsec:formal}
Let $\mathcal{X}$ be a locally noetherian formal scheme over $\mathcal{S}$
and $\mathcal{I}=\mathcal{I}_{\mathcal{X}}$ the largest ideal of definition of $\mathcal{X}$.
We put $\mathcal{X}_\red=(\mathcal{X},\mathcal{O}_{\mathcal{X}}/\mathcal{I})$, which is a locally noetherian
reduced scheme. The following lemma might be well-known:

\begin{lem}\label{lem:X-affine}
 If $\mathcal{X}_\red$ is affine, the formal scheme $\mathcal{X}$ is also affine.
\end{lem}

\begin{prf}
 Denote by $X_n$ the scheme $(\mathcal{X},\mathcal{O}_{\mathcal{X}}/\mathcal{I}^{n+1})$.
 Then, by \cite[I, (10.6.2)]{EGA},
 $\mathcal{X}=\varinjlim_n X_n$ (inductive limit in the category of formal schemes).
 We know that $X_0=\mathcal{X}_\red$ is affine. Therefore, by \cite[(2.3.5)]{EGA1-new},
 $X_n$ is an affine scheme
 for every $n\ge 0$. Put $A_n=\Gamma(\mathcal{X},\mathcal{O}_{\mathcal{X}}/\mathcal{I}^{n+1})$
 and $A=\varprojlim_n A_n$.
 Then $A$ is an admissible ring with a fundamental system of ideals of definitions 
 $\{\Ker (A\longrightarrow A_n)\}$, and we have $\mathcal{X}=\Spf A$.
\end{prf}

\begin{rem}\label{rem:A-I-adic}
 In the proof of the lemma above, put $I=\Ker (A\longrightarrow A_0)$. Then we have
 $I^{n+1}=\Ker (A\longrightarrow A_n)$ \cite[$\mathrm{0_I}$, (7.2.7)]{EGA}.
 In particular, the topology of $A$ coincides with the $I$-adic topology.
\end{rem}

\begin{defn}
 For an open subscheme $U$ of $\mathcal{X}_\red$, we denote by $\mathcal{X}_{/U}$ the open formal subscheme of
 $\mathcal{X}$ whose underlying space is $U$. If $U$ is affine, 
 $\mathcal{X}_{/U}$ is also affine by Lemma \ref{lem:X-affine}.
 Then we put $A_U=\Gamma(\mathcal{X}_{/U},\mathcal{O}_{\mathcal{X}})$,
 $I_U=\Gamma(\mathcal{X}_{/U},\mathcal{I})$
 and $\widehat{X}_{\!/U}=\Spec A_U$. Since $A_U$ is an $R$-algebra, $\widehat{X}_{\!/U}$ has a natural structure
 of an $S$-scheme.

 More generally, let $U$ be a (possibly infinite) disjoint union of affine open subschemes of $\mathcal{X}_\red$
 and $U=\coprod_{i\in I}U_i$ the decomposition into the connected components. Then we put
 $\widehat{X}_{\!/U}=\coprod_{i\in I}\widehat{X}_{\!/U_i}$. This extends the construction above,
 since $U\longmapsto \widehat{X}_{\!/U}$ is compatible with disjoint union.
\end{defn}

% \begin{defn}\label{defn:X_U}
%  For an \'etale morphism $U\longrightarrow \mathcal{X}_\red$, $\mathcal{X}_{/U}$ denotes the unique formal scheme
%  over $\mathcal{X}$ such that
%  \begin{itemize}
%   \item $\mathcal{X}_{/U}\longrightarrow \mathcal{X}$ is \'etale,
%   \item and $\mathcal{X}_{/U}\times_{\mathcal{X}}\mathcal{X}_\red\longrightarrow \mathcal{X}_\red$ coincides
% 	with $U\longrightarrow \mathcal{X}_\red$.
%  \end{itemize}
%  Note that $(\mathcal{X}_{/U})_\red=U$.
% 
%  If moreover $U$ is affine, we put $A_U=\Gamma(\mathcal{X}_{/U},\mathcal{O}_{\mathcal{X}_{/U}})$,
%  $I_U=\Gamma(\mathcal{X}_{/U},\mathcal{I}_{\mathcal{X}_{/U}})$
%  and $\widehat{X}_{\!/U}=\Spec A_U$. Since $A_U$ is an $R$-algebra, $\widehat{X}_{\!/U}$ has a natural structure
%  of an $S$-scheme.
% \end{defn}

% \begin{exa}
%  If $U$ is an open subscheme of $\mathcal{X}_\red$, then $\mathcal{X}_{/U}$ is the open subscheme of
%  $\mathcal{X}$ whose underlying space is $U$. More generally, let $U=\coprod_{i\in I}U_i$ be
%  a (possibly infinite) disjoint union of affine open subschemes of $\mathcal{X}_\red$,
%  then $\mathcal{X}_{/U}=\coprod_{i\in I}\mathcal{X}_{/U_i}$, where $\mathcal{X}_{/U_i}$ can be described
%  as above.
%  In fact, these examples are all what we need.
% \end{exa}

\begin{lem}\label{lem:X_U-funct}
 Let $f\colon \mathcal{Y}\longrightarrow \mathcal{X}$ be a morphism between locally noetherian formal schemes
 over $\mathcal{S}$,
 $U$ (resp.\ $V$) a disjoint union of affine open subschemes of $\mathcal{X}_\red$ (resp.\ $\mathcal{Y}_\red$)
 and $f'\colon V\longrightarrow U$ a morphism of schemes which makes the following diagram
 commutative:
 \[
  \xymatrix{%
 V\ar[r]^-{f'}\ar[d]& U\ar[d]\\ \mathcal{Y}_\red\ar[r]^-{f_\red}& \mathcal{X}_\red\lefteqn{.}
 }
 \]
 Then we have a natural $S$-morphism $\widehat{Y}_{\!/V}\longrightarrow \widehat{X}_{\!/U}$.
\end{lem}

\begin{prf}
 Let $U=\coprod_{i\in I}U_i$ and $V=\coprod_{j\in J}V_j$
 be the decompositions into the connected components. For every $j\in J$, there exists unique $i\in I$
 such that $f'(V_j)\subset U_i$. Therefore $f'$ induces the morphism of affine formal schemes
 $\Spf B_{V_j}=\mathcal{Y}_{/V_j}\longrightarrow \mathcal{X}_{/U_i}=\Spf A_{U_i}$.
 This induces the continuous $R$-algebra homomorphism $A_{U_i}\longrightarrow B_{V_j}$ 
 and the morphism of $S$-schemes $\widehat{Y}_{\!/V_j}=\Spec B_{V_j}\longrightarrow \Spec A_{U_i}=\widehat{X}_{\!/U_i}\hooklongrightarrow \widehat{X}_{\!/U}$. This gives the morphism of $S$-schemes
 $\widehat{Y}_{\!/V}=\coprod_{j\in J}\widehat{Y}_{\!/V_j}\longrightarrow \widehat{X}_{\!/U}$.
\end{prf}

\begin{lem}\label{lem:X_U-flat}
 Let $U$ be an affine open subscheme of $\mathcal{X}_{\red}$.
 \begin{enumerate}
  \item For an affine open subscheme $U'$ of $U$, 
	the natural morphism $\widehat{X}_{\!/U'}\longrightarrow \widehat{X}_{\!/U}$ is flat.
  \item For an affine open covering $\{U_i\}_{i\in I}$ of $U$, 
	the morphism $\coprod_{i\in I}\widehat{X}_{\!/U_i}\longrightarrow \widehat{X}_{\!/U}$ is
	faithfully flat.
 \end{enumerate}
\end{lem}

\begin{prf}
 \begin{enumerate}
  \item By Remark \ref{rem:A-I-adic}, we have $A_U=\varprojlim_{n}A_U/I_U^n$
	and $A_{U'}=\varprojlim_{n}A_{U'}/I_U^nA_{U'}$.
	On the other hand, $\Spec A_{U'}/I_U^nA_{U'}\longrightarrow \Spec A_U/I_U^n$ is flat since
	it is an open immersion.
	Therefore, by the local criterion of flatness 
	\cite[$\mathrm{0_{III}}$, (10.2.2)]{EGA}, $\Spec A_{U'}\longrightarrow \Spec A_U$ is flat, as desired.
  \item By i), the morphism is flat. On the other hand, the image contains the closed
	subset $U\subset \widehat{X}_{\!/U}$. Since every closed point of $\widehat{X}_{\!/U}$ lies in $U$,
	we have the surjectivity of $\coprod_{i\in I}\widehat{X}_{\!/U_i}\longrightarrow \widehat{X}_{\!/U}$.
 \end{enumerate}
\end{prf}

Next we introduce the notion of excellent formal schemes. 

\begin{defn}
 We say that $\mathcal{X}$ is excellent if for every affine open subscheme $U$ of $\mathcal{X}_\red$
 the ring $A_U$ is excellent.
\end{defn}

The following proposition gives a lot of examples of excellent formal schemes:

\begin{prop}\label{prop:special-exc}
 Every special formal scheme over $\mathcal{S}$ in the sense of Berkovich (\cf \cite{MR1395723}) is excellent.
 In particular, every formal scheme obtained by completing an $S$-scheme $X$ locally of finite type 
 along a closed subscheme of $X_s$ is excellent.
\end{prop}

\begin{prf}
 Let $\mathcal{X}$ be a special formal scheme over $\mathcal{S}$. Then, by definition, there exists an affine open covering $\{U_i\}$ of $\mathcal{X}_\red$
 such that $A_{U_i}$ is a special $R$-algebra for every $i$. We will observe that for every affine open subscheme $U$ of $\mathcal{X}_\red$,
 the ring $A_U$ is a special $R$-algebra. By \cite[Lemma 1.2]{MR1395723}, it suffices to show that $A_U/I_U^2$ is a finitely generated $R$-algebra. Take an affine open covering $\{V_j\}$ of $U$ such that for each $V_j$
 is contained in $U_i$ for some $i$.
 Since $\Spec A_U/I_U^2$ is covered by open subschemes $\Spec A_{V_j}/I_{V_j}^2$,
 it is sufficient to prove that $A_{V_j}/I_{V_j}^2$ is
 a finitely generated $R$-algebra. This is clear, since we have a quasi-compact open immersion 
 $\Spec A_{V_j}/I_{V_j}^2\hooklongrightarrow\Spec A_{U_i}/I_{U_i}^2$ for some $i$.

 Therefore, by \cite[Lemma 1.2]{MR1395723}, $A_U$ is a quotient of $R\langle T_1,\ldots,T_m\rangle[[S_1,\ldots,S_n]]$ for some integers $m$ and $n$.
 By \cite[Theorem 9]{MR0407007} (where $R$ has mixed characteristic) and \cite[Proposition 7]{MR0376677} (where $R$ has equal characteristic),
 the $R$-algebra $R\langle T_1,\ldots,T_m\rangle[[S_1,\ldots,S_n]]$ is excellent. Thus $A_U$ is also excellent.
\end{prf}

\begin{prop}\label{prop:regular}
 Assume that $\mathcal{X}$ is excellent. Let $U$ and $U'$ be disjoint unions of affine open subschemes of
 $\mathcal{X}_\red$ and $U'\longrightarrow U$ a morphism over $\mathcal{X}_\red$.
 Then the natural morphism $\widehat{X}_{\!/U'}\longrightarrow \widehat{X}_{\!/U}$ is regular.
\end{prop}

\begin{prf}
 In the same way as in the proof of Lemma \ref{lem:X_U-funct}, we can reduce to the case
 where $U$ and $U'$ are affine open subschemes of $\mathcal{X}_\red$ with $U'\subset U$.

 First we consider the case where $U'=\mathfrak{D}(f)_\red$ with $f\in A_U$. Then $A_{U'}$ is the $I_U$-adic completion of $(A_U)_f$.
 Since $A_U$ is excellent, so is $(A_U)_f$, and we have the regularity of the morphism 
 $\widehat{X}_{\!/U'}=\Spec A_{U'}\longrightarrow \Spec (A_U)_f$ \cite[IV, (7.8.3) (v)]{EGA}.
 On the other hand, it is clear that the morphism $\Spec (A_U)_f\longrightarrow \Spec A_U=\widehat{X}_{\!/U}$ is regular.
 This completes the proof for the case $U'=\mathfrak{D}(f)_\red$.

 In the general case, we can cover $U'$ by open subschemes of the form $\mathfrak{D}(f_i)_\red$ with $f_i\in A_U$.
 Then $\coprod_{i}\widehat{X}_{\!/\mathfrak{D}(f_i)_\red}\longrightarrow \widehat{X}_{\!/U}$ is regular.
 On the other hand, by Lemma \ref{lem:X_U-flat} ii), the morphism $\coprod_{i}\widehat{X}_{\!/\mathfrak{D}(f_i)_\red}\longrightarrow \widehat{X}_{\!/U'}$ is
 faithfully flat. Therefore, $\widehat{X}_{\!/U'}\longrightarrow \widehat{X}_{\!/U}$ is regular by
 \cite[IV, (7.3.8) (iv$'$)]{EGA}.
\end{prf}

Let $\mathcal{Z}$ be a closed subscheme of $\mathcal{X}$ and $U$ a disjoint union of affine open subschemes of
$\mathcal{X}_\red$. Then, as $\mathcal{Z}_\red\times_{\mathcal{X}_\red} U$ is a disjoint union of
affine open subschemes of $\mathcal{Z}_\red$, we can consider a formal scheme 
$\mathcal{Z}_{/\mathcal{Z}_\red\times_{\mathcal{X}_\red} U}$ and
a scheme $\widehat{Z}_{/\mathcal{Z}_\red\times_{\mathcal{X}_\red} U}$. For simplicity, we denote them
by $\mathcal{Z}_{/U}$ and $\widehat{Z}_{\!/U}$, respectively.
The former is a closed formal subscheme of $\mathcal{X}_{/U}$ and the latter is
a closed subscheme of $\widehat{X}_{\!/U}$.

\begin{lem}\label{lem:closed-cart}
 Let $f\colon \mathcal{Y}\longrightarrow \mathcal{X}$ and $U$, $V$ be as in Lemma \ref{lem:X_U-funct}.
 Put $\mathcal{Z}'=\mathcal{Y}\times_{\mathcal{X}}\mathcal{Z}$, which is a closed subscheme of $\mathcal{Y}$.
 Then the following diagrams are cartesian:
 \[
 \xymatrix{%
 \mathcal{Z}'_{/V}\ar[r]\ar[d]& \mathcal{Y}_{/V}\ar[d]\\
 \mathcal{Z}_{/U}\ar[r]& \mathcal{X}_{/U}\lefteqn{,}
 }\qquad\qquad
 \xymatrix{%
 \widehat{Z}'_{\!/V}\ar[r]\ar[d]& \widehat{Y}_{\!/V}\ar[d]\\
 \widehat{Z}_{\!/U}\ar[r]& \widehat{X}_{\!/U}\lefteqn{.}
 }
 \]
\end{lem}

\begin{prf}
 We may assume that $U$ is an affine open subscheme of $\mathcal{X}_\red$
 and $V$ is an affine open subscheme of $\mathcal{Y}_\red$.  
 It is clear that the left diagram is cartesian. Let $J$ be the ideal of $A_U$ such that
 $\mathcal{Z}_{/U}=\Spf A_U/J$. Then $\mathcal{Z}'_{/V}=\Spf A_U/J\times_{\Spf A_U}\Spf B_V=\Spf B_V/JB_V$,
 where $B_V=\Gamma(\mathcal{Y}_{/V},\mathcal{O}_{\mathcal{Y}})$; note that $B_V/JB_V$ is complete since
 $B_V$ is noetherian. Thus $\widehat{Z}'_{\!/V}=\Spec B_V/JB_V$ and the right diagram is also cartesian.
\end{prf}

\begin{defn}
 Let $\mathcal{Z}=(\mathcal{Z}_1,\mathcal{Z}_2)$ be a pair of closed formal subschemes of $\mathcal{X}$
 such that $\mathcal{Z}_1\subset \mathcal{Z}_2$. 
 Then $(\mathcal{Z}_1)_{/U}$ is a closed formal subscheme of $(\mathcal{Z}_2)_{/U}$ and
 $(\widehat{Z}_1)_{/U}$ is a closed subscheme of $(\widehat{Z}_2)_{/U}$.  
 We put $\widehat{Z}_{\!/U}=(\widehat{Z}_2)_{/U}\setminus (\widehat{Z}_1)_{/U}$,
 which is a locally closed subscheme of $\widehat{X}_{\!/U}$.
\end{defn}

Under the situation of Lemma \ref{lem:X_U-funct}, put $\mathcal{Z}'=(\mathcal{Y}\times_{\mathcal{X}}\mathcal{Z}_1,\mathcal{Y}\times_{\mathcal{X}}\mathcal{Z}_2)$. Then the natural $S$-morphism
$\widehat{Y}_{\!/V}\longrightarrow \widehat{X}_{\!/U}$ induces a morphism
$\widehat{Z}'_{\!/V}\longrightarrow \widehat{Z}_{\!/U}$; indeed, by Lemma \ref{lem:closed-cart},
the following diagram is cartesian:
\[
  \xymatrix{%
 (\widehat{Z}'_1)_{/V}\ar[r]\ar[d]& (\widehat{Z}'_2)_{/V}\ar[d]\\
 (\widehat{Z}_1)_{/U}\ar[r]& (\widehat{Z}_2)_{/U}\lefteqn{.}
 }
\]

\section{Variant of formal nearby cycles $R\Psi_{\!\mathcal{X},\mathcal{Z}}\Lambda$}\label{sec:vfnc}
\subsection{Construction}\label{subsec:construction}
As in the previous section, let $\mathcal{X}$ be a locally noetherian formal scheme over $\mathcal{S}$.
Moreover assume that $\mathcal{X}$ is excellent and $\mathcal{X}_\red$ is separated.
Let $\mathcal{Z}=(\mathcal{Z}_1,\mathcal{Z}_2)$ be a pair of closed formal subschemes of $\mathcal{X}_s$
with $\mathcal{Z}_1\subset \mathcal{Z}_2$.
In this subsection, we will construct the object $R\Psi_{\!\mathcal{X},\mathcal{Z}}\Lambda$ of
$D^+(\mathcal{X}_\red,\Lambda)$.

Let $U$ be an affine open covering of $\mathcal{X}_\red$
and $a\colon U_\bullet\longrightarrow \mathcal{X}_\red$ the associated hypercovering
(\cf Section \ref{subsec:open-cov}). Since $\mathcal{X}_\red$ is separated, $U_m$ is a disjoint union
of affine open subschemes of $\mathcal{X}_\red$ for every $m\ge 0$.
Therefore, we have an $S$-scheme $\widehat{X}_{\!/U_m}$ for every $m$. By Lemma \ref{lem:X_U-funct}, 
$\widehat{X}_{\!/U_\bullet}$ naturally have a structure of a simplicial $S$-scheme.
Similarly, we have the simplicial $s$-scheme $\widehat{Z}_{\!/U_\bullet}$. Denote the natural
immersion $U_\bullet\longrightarrow (\widehat{X}_{\!/U_\bullet})_s$ 
(resp.\ $\widehat{Z}_{\!/U_\bullet}\longrightarrow (\widehat{X}_{\!/U_\bullet})_s$) by $i$ (resp.\ $j$)
or by $i_U$ (resp.\ $j_U$) if we need to indicate $U$.
\[
 \xymatrix{%
 U_\bullet\ar[r]^-{i}\ar[d]^-{a}& (\widehat{X}_{\!/U_\bullet})_s& \widehat{Z}_{\!/U_\bullet}\ar[l]_-{j}\\
 \mathcal{X}_\red
 }
\]

\begin{lem}\label{lem:Rpsi-cart}
 The complex $i^*Rj_*Rj^!R\psi_{\widehat{X}_{\!/U_\bullet}}\Lambda$ on $\mathcal{X}_\red$ lies
 in $D^+_{\mathrm{cart}}(\mathcal{X}_\red,\Lambda)$.
 Here $R\psi_{\widehat{X}_{\!/U_\bullet}}\colon D^+(\widehat{X}_{\!/U_\bullet},\Lambda)\longrightarrow D^+((\widehat{X}_{\!/U_\bullet})_s,\Lambda)$ is
 the simplicial version of the nearby cycle functor over $S$, whose definition is the obvious one.
\end{lem}

\begin{prf}
 Note that the restriction of $i^*Rj_*Rj^!R\psi_{\widehat{X}_{\!/U_\bullet}}\Lambda$ to $U_m$
 is $i_m^*Rj_{m*}Rj_m^!R\psi_{\widehat{X}_{\!/U_m}}\Lambda$. Therefore, 
 it suffices to show that, for every structure morphism $\phi\colon U_m\longrightarrow U_n$ of $U_\bullet$,
 the natural morphism $\phi^*i_n^*Rj_{n*}Rj_n^!R\psi_{\widehat{X}_{\!/U_n}}\Lambda\longrightarrow i_m^*Rj_{m*}Rj_m^!R\psi_{\widehat{X}_{\!/U_m}}\Lambda$ is an isomorphism.
 Consider the following commutative diagrams:
 \[
 \xymatrix{%
 U_m\ar[r]^-{i_m}\ar[d]^-{\phi}& (\widehat{X}_{\!/U_m})_s\ar[d]^-{\phi}\\
 U_n\ar[r]^-{i_n}& (\widehat{X}_{\!/U_n})_s\lefteqn{,}
 }
 \qquad\qquad
 \xymatrix{%
 \widehat{Z}_{\!/U_m}\ar[r]^-{j_m}\ar[d]^-{\phi}& (\widehat{X}_{\!/U_m})_s\ar[d]^-{\phi}\\
 \widehat{Z}_{\!/U_n}\ar[r]^-{j_n}& (\widehat{X}_{\!/U_n})_s\lefteqn{.}
 }
 \]
 The right diagram is cartesian by Lemma \ref{lem:closed-cart}.
 All morphisms $\phi$ in the diagrams above are regular by Proposition \ref{prop:regular}.
 Therefore,
 \begin{align*}
  \phi^*i_n^*Rj_{n*}Rj_n^!R\psi_{\widehat{X}_{\!/U_n}}\Lambda&=i_m^*\phi^*Rj_{n*}Rj_n^!R\psi_{\widehat{X}_{\!/U_n}}\Lambda\yrightarrow[(1)]{\cong}i_m^*Rj_{m*}\phi^*Rj_n^!R\psi_{\widehat{X}_{\!/U_n}}\Lambda\\
  &\yrightarrow[(2)]{\cong}i_m^*Rj_{m*}Rj_m^!\phi^*R\psi_{\widehat{X}_{\!/U_n}}\Lambda
  \yrightarrow[(3)]{\cong}i_m^*Rj_{m*}Rj_m^!R\psi_{\widehat{X}_{\!/U_m}}\Lambda,
 \end{align*}
 as desired. The isomorphy of (1) and (3) is a consequence of the regular base change theorem
 \cite[Corollary 7.1.6]{MR1360610}. The isomorphy of (2) follows from \cite[Corollaire 4.7]{Dualite}.
\end{prf}

\begin{defn}
 We define the object $R\Psi_{\!\mathcal{X},\mathcal{Z},U}\Lambda\in D^+(\mathcal{X}_\red,\Lambda)$ by
 \[
 R\Psi_{\!\mathcal{X},\mathcal{Z},U}\Lambda=Ra_*i^*Rj_*Rj^!R\psi_{\widehat{X}_{\!/U_\bullet}}\Lambda.
 \]
 By Lemma \ref{lem:Rpsi-cart}, it is characterized by the property 
 $a^*R\Psi_{\!\mathcal{X},\mathcal{Z},U}\Lambda=i^*Rj_*Rj^!R\psi_{\widehat{X}_{\!/U_\bullet}}\Lambda$. 
\end{defn}

We would like to observe that $R\Psi_{\!\mathcal{X},\mathcal{Z},U}\Lambda$ is independent of 
$U$ up to canonical isomorphism.
For this purpose, let $V$ be another affine open covering of $\mathcal{X}_\red$
and $b\colon V_\bullet\longrightarrow \mathcal{X}_\red$ the associated hypercovering.
We will construct an isomorphism 
$\lambda_{VU}\colon R\Psi_{\!\mathcal{X},\mathcal{Z},U}\Lambda\longrightarrow R\Psi_{\!\mathcal{X},\mathcal{Z},V}\Lambda$ by means of the bisimplicial technique.
Let $c\colon W_{\bullet\bullet}\longrightarrow \mathcal{X}_\red$ be the augmented bisimplicial scheme
associated with $a\colon U_\bullet\longrightarrow \mathcal{X}_\red$ and
$b\colon V_\bullet\longrightarrow \mathcal{X}_\red$ (\cf Section \ref{subsec:open-cov}).
For every $m,n\ge 0$, $W_{mn}=U_m\times_{\mathcal{X}_\red}V_n$ is a disjoint union
of affine open subschemes of $\mathcal{X}_\red$, since $\mathcal{X}_\red$ is separated.
We have the first projection $p\colon W_{\bullet\bullet}\longrightarrow U_\bullet$
and the second projection $q\colon W_{\bullet\bullet}\longrightarrow V_\bullet$.

In the same way as $\widehat{X}_{\!/U_\bullet}$, we can construct the bisimplicial $S$-schemes
$\widehat{X}_{\!/W_{\bullet\bullet}}$ and $\widehat{Z}_{\!/W_{\bullet\bullet}}$.
We have natural morphisms between \'etale sites illustrated in the following diagrams:
\[
 \xymatrix{%
 W_{\bullet\bullet}\ar[r]^-{i_W}\ar[d]^-{p}\ar@/_20pt/[dd]_-{c}& (\widehat{X}_{\!/W_{\bullet\bullet}})_s\ar[d]^-{p}& \widehat{Z}_{\!/W_{\bullet\bullet}}\ar[l]_-{j_W}\ar[d]^-{p}\\
 U_\bullet\ar[r]^-{i_U}\ar[d]^-{a}& (\widehat{X}_{\!/U_\bullet})_s& \widehat{Z}_{\!/U_\bullet}\ar[l]_-{j_U}\\
 \mathcal{X}_\red\lefteqn{,}
 }
 \qquad
  \xymatrix{%
 W_{\bullet\bullet}\ar[r]^-{i_W}\ar[d]^-{q}\ar@/_20pt/[dd]_-{c}& (\widehat{X}_{\!/W_{\bullet\bullet}})_s\ar[d]^-{q}& \widehat{Z}_{\!/W_{\bullet\bullet}}\ar[l]_-{j_W}\ar[d]^-{q}\\
 V_\bullet\ar[r]^-{i_V}\ar[d]^-{b}& (\widehat{X}_{\!/V_\bullet})_s& \widehat{Z}_{\!/V_\bullet}\ar[l]_-{j_V}\\
 \mathcal{X}_\red\lefteqn{.}
 }
\]
We can construct the object $i_W^*Rj_{W*}Rj^!_WR\psi_{\widehat{X}_{\!/W_{\bullet\bullet}}}\Lambda$ of
$D^+(W_{\bullet\bullet},\Lambda)$ and the natural morphisms
\begin{align*}
 &p^*i_U^*Rj_{U*}Rj^!_UR\psi_{\widehat{X}_{\!/U_{\bullet}}}\Lambda
 \yrightarrow{\varphi_p} i_W^*Rj_{W*}Rj^!_WR\psi_{\widehat{X}_{\!/W_{\bullet\bullet}}}\Lambda,\\
 &q^*i_V^*Rj_{V*}Rj^!_VR\psi_{\widehat{X}_{\!/V_{\bullet}}}\Lambda
 \yrightarrow{\varphi_q} i_W^*Rj_{W*}Rj^!_WR\psi_{\widehat{X}_{\!/W_{\bullet\bullet}}}\Lambda.
\end{align*}

\begin{lem}\label{lem:bisimplicial}
 The two morphisms $\varphi_p$, $\varphi_q$ are isomorphisms.
\end{lem}

\begin{prf}
 It suffices to show that $\varphi_p\vert_{W_{mn}}$ is an isomorphism for every $m,n\ge 0$.
 We have 
 \begin{align*}
  (p^*i_U^*Rj_{U*}Rj^!_UR\psi_{\widehat{X}_{\!/U_{\bullet}}}\Lambda)\vert_{W_{mn}}&=p_{mn}^*i_{U,m}^*Rj_{U,m*}Rj_{U,m}^!R\psi_{\widehat{X}_{\!/U_m}}\Lambda,\\
  (i_W^*Rj_{W*}Rj^!_WR\psi_{\widehat{X}_{\!/W_{\bullet\bullet}}}\Lambda)\vert_{W_{mn}}&=i_{W,mn}^*Rj_{W,mn*}Rj_{W,mn}^!R\psi_{\widehat{X}_{\!/W_{mn}}}\Lambda,
 \end{align*}
 and $\varphi_p\vert_{W_{mn}}$ can be identified with the morphism naturally induced from the diagram below,
 whose right rectangle is cartesian:
 \[
 \xymatrix{%
 W_{mn}\ar[r]^-{i_{W,mn}}\ar[d]^-{p_{mn}}& (\widehat{X}_{\!/W_{mn}})_s\ar[d]^-{p_{mn}}& \widehat{Z}_{\!/W_{mn}}\ar[l]_-{j_{W,mn}}\ar[d]^-{p_{mn}}\\
 U_m\ar[r]^-{i_{U,m}}& (\widehat{X}_{\!/U_m})_s& \widehat{Z}_{\!/U_m}\ar[l]_-{j_{U,m}}\lefteqn{.}
 }
 \]
 Hence we can show the isomorphy of $\varphi_p\vert_{W_{mn}}$ exactly in the same way as in the proof of
 Lemma \ref{lem:Rpsi-cart} by using the regularity of 
 $p_{mn}\colon \widehat{X}_{\!/W_{mn}}\longrightarrow \widehat{X}_{\!/U_m}$
 (Proposition \ref{prop:regular}).
\end{prf}

By this lemma, we have the isomorphism
\[
 c^*R\Psi_{\!\mathcal{X},\mathcal{Z},U}\Lambda=p^*i_U^*Rj_{U*}Rj_U^!R\psi_{\widehat{X}_{\!/U_\bullet}}\Lambda
 \yrightarrow{\varphi_q^{-1}\circ \varphi_p}
 q^*i_V^*Rj_{V*}Rj^!_VR\psi_{\widehat{X}_{\!/V_{\bullet}}}\Lambda
 =c^*R\Psi_{\!\mathcal{X},\mathcal{Z},V}\Lambda.
\]
Since $c^*$ is fully faithful (Lemma \ref{lem:bisimplicial-descent}), there exists a unique
isomorphism 
\[
 \lambda_{VU}\colon R\Psi_{\!\mathcal{X},\mathcal{Z},U}\Lambda\yrightarrow{\cong} R\Psi_{\!\mathcal{X},\mathcal{Z},V}\Lambda
\]
such that $c^*(\lambda_{VU})=\varphi_q^{-1}\circ \varphi_p$.

\begin{lem}\label{lem:inductive}
 \begin{enumerate}
  \item For an affine open covering $U$ of $\mathcal{X}_\red$, we have $\lambda_{UU}=\id$.
  \item For affine open coverings $U$, $V$, $W$ of $\mathcal{X}_\red$, 
	we have $\lambda_{WV}\circ \lambda_{VU}=\lambda_{WU}$.
 \end{enumerate}
\end{lem}

\begin{prf}
By using the trisimplicial scheme
$U_\bullet\times_{\mathcal{X}_\red}V_\bullet\times_{\mathcal{X}_\red}W_\bullet$,
it is straightforward to prove ii).
 For i), note that $\lambda_{UU}\circ \lambda_{UU}=\lambda_{UU}$ implies $\lambda_{UU}=\id$,
 since $\lambda_{UU}$ is an isomorphism.
\end{prf}

\begin{defn}
 We put $R\Psi_{\!\mathcal{X},\mathcal{Z}}\Lambda=\varinjlim_{U}R\Psi_{\!\mathcal{X},\mathcal{Z},U}\Lambda\in D^+(\mathcal{X}_\red,\Lambda)$. Here the inductive limit is taken over the small category as follows:
 \begin{itemize}
  \item the objects are affine open coverings of $\mathcal{X}_\red$,
  \item and $\Hom(U,V)$ consists of one element for every objects $U$, $V$.
 \end{itemize}
 Lemma \ref{lem:inductive} ensures that $(R\Psi_{\!\mathcal{X},\mathcal{Z},U}\Lambda)_U$ forms
 an inductive system.
 Since $\lambda_{VU}$ is an isomorphism, the existence of the inductive limit is immediate.
 Note that the canonical morphism $\lambda_U\colon R\Psi_{\!\mathcal{X},\mathcal{Z},U}\Lambda\longrightarrow R\Psi_{\!\mathcal{X},\mathcal{Z}}\Lambda$ is an isomorphism.

 If $\mathcal{Z}=(\varnothing, \mathcal{X}_s)$ (resp.\ $\mathcal{Z}=(\varnothing, \mathcal{X}_\red)$),
 then $R\Psi_{\!\mathcal{X},\mathcal{Z}}\Lambda$ is denoted by $R\Psi_{\!\mathcal{X}}\Lambda$
 (resp.\ $R\Psi_{\!\mathcal{X},c}\Lambda$). 
\end{defn}

\subsection{Basic properties}\label{sec:properties}
In this subsection, we gather some basic properties of $R\Psi_{\!\mathcal{X},\mathcal{Z}}\Lambda$.
Let $\mathcal{X}$ be the same as in the previous section.

\subsubsection{Excision triangle}
\begin{prop}\label{prop:excision}
 Let $\mathcal{Z}_1$, $\mathcal{Z}_2$ and $\mathcal{Z}_3$ be closed formal subschemes of $\mathcal{X}_s$
 with $\mathcal{Z}_1\subset \mathcal{Z}_2\subset\mathcal{Z}_3$, and
 put $\mathcal{Z}=(\mathcal{Z}_1,\mathcal{Z}_2)$, $\mathcal{Z}'=(\mathcal{Z}_1,\mathcal{Z}_3)$ and
 $\mathcal{Z}''=(\mathcal{Z}_2,\mathcal{Z}_3)$. Then we have a natural distinguished triangle
 \[
  R\Psi_{\!\mathcal{X},\mathcal{Z}}\Lambda\longrightarrow R\Psi_{\!\mathcal{X},\mathcal{Z}'}\Lambda
 \longrightarrow R\Psi_{\!\mathcal{X},\mathcal{Z}''}\Lambda\longrightarrow R\Psi_{\!\mathcal{X},\mathcal{Z}}\Lambda[1].
 \]
\end{prop}

\begin{prf}
 Take an affine open covering $U$ of $\mathcal{X}_{\red}$ and 
 let $a\colon U_\bullet\longrightarrow \mathcal{X}_\red$ be the associated hypercovering.
 Then we have three locally closed subschemes of $(\widehat{X}_{\!/U_\bullet})_s$:
 \[
  j\colon \widehat{Z}_{\!/U_\bullet}\hooklongrightarrow (\widehat{X}_{\!/U_\bullet})_s,
 \quad j'\colon \widehat{Z}'_{\!/U_\bullet}\hooklongrightarrow (\widehat{X}_{\!/U_\bullet})_s,
 \quad j''\colon \widehat{Z}''_{\!/U_\bullet}\hooklongrightarrow (\widehat{X}_{\!/U_\bullet})_s.
 \]
 By the definition, $\widehat{Z}_{\!/U_\bullet}$ is a closed simplicial subscheme of
 $\widehat{Z}'_{\!/U_\bullet}$ whose complement coincides with $\widehat{Z}''_{\!/U_\bullet}$.
 Thus we have the following distinguished triangle:
 \[
  Rj_*Rj^!R\psi_{\widehat{X}_{\!/U_\bullet}}\Lambda\longrightarrow 
 Rj'_*Rj'^!R\psi_{\widehat{X}_{\!/U_\bullet}}\Lambda\longrightarrow 
 Rj''_*Rj''^!R\psi_{\widehat{X}_{\!/U_\bullet}}\Lambda\longrightarrow 
 Rj_*Rj^!R\psi_{\widehat{X}_{\!/U_\bullet}}\Lambda[1].
 \]
 By taking $Ra_*i^*$, we get the desired distinguished triangle.
\end{prf}

\subsubsection{Functoriality}
Let $\mathcal{X}'$ be another locally noetherian excellent formal scheme over $\mathcal{S}$
such that $\mathcal{X}'_\red$ is separated, and
$f\colon \mathcal{X}'\longrightarrow \mathcal{X}$ a morphism over $\mathcal{S}$.
For a pair $\mathcal{Z}=(\mathcal{Z}_1,\mathcal{Z}_2)$ of closed formal subschemes of $\mathcal{X}_s$
with $\mathcal{Z}_1\subset \mathcal{Z}_2$, put $\mathcal{Z}'_1=\mathcal{X}'\times_{\mathcal{X}}\mathcal{Z}_1$, 
$\mathcal{Z}'_2=\mathcal{X}'\times_{\mathcal{X}}\mathcal{Z}_2$ and
$\mathcal{Z}'=(\mathcal{Z}'_1, \mathcal{Z}'_2)$.

\begin{prop}\label{prop:functoriality}
 Under the setting above, we may construct the natural morphism
 $\varphi_f\colon f_\red^*R\Psi_{\!\mathcal{X},\mathcal{Z}}\Lambda\longrightarrow R\Psi_{\!\mathcal{X}',\mathcal{Z}'}\Lambda$. This is compatible with composition.
\end{prop}

\begin{prf}
 We can take an affine open covering $U=\{U_i\}_{i\in I}$ of $\mathcal{X}_\red$,
 an affine open covering $U'=\{U'_{i'}\}_{i'\in I'}$ of $\mathcal{X}'_\red$ and 
 a map $\alpha\colon I'\longrightarrow I$ such that $f_\red(U'_{i'})\subset U_{\alpha(i')}$
 for every $i'\in I'$. Let $a\colon U_\bullet\longrightarrow \mathcal{X}_\red$
 (resp.\ $a'\colon U'_\bullet\longrightarrow \mathcal{X}'_\red$) be the hypercovering associated with $U$
 (resp.\ $U'$). Then $\alpha$ induces the morphism of simplicial schemes 
 $f_\alpha\colon U'_\bullet\longrightarrow U_\bullet$ that makes the following diagram commutative:
 \[
  \xymatrix{%
 U'_\bullet\ar[r]^-{f_\alpha}\ar[d]^-{a'}& U_\bullet\ar[d]^-{a}\\
 \mathcal{X}'_\red\ar[r]^-{f_\red}& \mathcal{X}_\red\lefteqn{.}
 }
 \]
 Therefore, by Lemma \ref{lem:X_U-funct}, we have the morphism of simplicial schemes
 $f_\alpha\colon \widehat{X}'_{\!/U'_\bullet}\longrightarrow \widehat{X}_{\!/U_\bullet}$
 and the following commutative diagram whose right rectangle is cartesian:
 \[
 \xymatrix{%
 U'_\bullet\ar[r]^-{i'}\ar[d]^-{f_\alpha}& (\widehat{X}'_{\!/U'_\bullet})_s\ar[d]^-{f_\alpha}& \widehat{Z}'_{\!/U'_\bullet}\ar[l]_-{j'}\ar[d]^-{f_\alpha}\\
 U_\bullet\ar[r]^-{i}& (\widehat{X}_{\!/U_\bullet})_s& \widehat{Z}_{\!/U_\bullet}\ar[l]_-{j}\lefteqn{.}
 }
 \]
 Thus we have natural morphisms
 \begin{align*}
  f_\alpha^*i^*Rj_*Rj^!R\psi_{\widehat{X}_{\!/U_\bullet}}\Lambda&=i'^*f_\alpha^*Rj_*Rj^!R\psi_{\widehat{X}_{\!/U_\bullet}}\Lambda\longrightarrow i'^*Rj'_*f_\alpha^*Rj^!R\psi_{\widehat{X}_{\!/U_\bullet}}\Lambda\\
  &\longrightarrow i'^*Rj'_*Rj'^!f_\alpha^*R\psi_{\widehat{X}_{\!/U_\bullet}}\Lambda
  \longrightarrow i'^*Rj'_*Rj'^!R\psi_{\widehat{X}'_{\!/U'_\bullet}}\Lambda.
 \end{align*}
 We will denote it by $\widetilde{\varphi}_\alpha$. This induces the morphism
 \begin{align*}
  a'^*f_\red^*R\Psi_{\!\mathcal{X},\mathcal{Z}}\Lambda&=f_\alpha^*a^*R\Psi_{\!\mathcal{X},\mathcal{Z}}\Lambda
  \yleftarrow[\cong]{f_\alpha^*a^*\lambda_U}f_\alpha^*a^*R\Psi_{\!\mathcal{X},\mathcal{Z},U}\Lambda=f_\alpha^*i^*Rj_*Rj^!R\psi_{\widehat{X}_{\!/U_\bullet}}\Lambda\\
  &\yrightarrow{\widetilde{\varphi}_\alpha} i'^*Rj'_*Rj'^!R\psi_{\widehat{X}'_{\!/U'_\bullet}}\Lambda
  =a'^*R\Psi_{\!\mathcal{X}',\mathcal{Z}',U'}\Lambda\yrightarrow[\cong]{a'^*\lambda'_{U'}}a'^*R\Psi_{\!\mathcal{X}',\mathcal{Z}'}\Lambda.
 \end{align*}
 We define $\varphi_\alpha$ as the unique morphism such that $a'^*(\varphi_\alpha)$ coincides with the morphism above.

 We will prove that this is independent of the choice of $(U,U',\alpha)$. Let $(V,V',\beta)$ be another triple
 and denote the hypercovering associated with $V$ (resp.\ $V'$) by $b\colon V\longrightarrow \mathcal{X}_\red$
 (resp.\ $b'\colon V'\longrightarrow \mathcal{X}'_\red$). Put $W_{\bullet\bullet}=U_\bullet\times_{\mathcal{X}_\red}V_\bullet$, $W'_{\bullet\bullet}=U'_\bullet\times_{\mathcal{X}'_\red}V'_\bullet$ and
 let $c\colon W_{\bullet\bullet}\longrightarrow \mathcal{X}_\red$ and
 $c'\colon W'_{\bullet\bullet}\longrightarrow \mathcal{X}'_\red$ be the natural augmentations.
 Then, as in Section \ref{subsec:construction}, we have the following diagrams:
 \[
  \xymatrix{%
 U'_\bullet\ar[d]^-{f_\alpha}& W'_{\bullet\bullet}\ar[l]_-{p'}\ar[d]^-{f_{\alpha\beta}}\ar[r]^-{q'}&
 V'_\bullet\ar[d]^-{f_\beta}\\
 U_\bullet& W_{\bullet\bullet}\ar[l]_-{p}\ar[r]^-{q'}&
 V_\bullet\lefteqn{,}
 }\quad
 \xymatrix{%
 \widehat{X}'_{\!U'_\bullet}\ar[d]^-{f_\alpha}& \widehat{X}'_{\!W'_{\bullet\bullet}}\ar[l]_-{p'}\ar[d]^-{f_{\alpha\beta}}\ar[r]^-{q'}&
 \widehat{X}'_{\!V'_\bullet}\ar[d]^-{f_\beta}\\
 \widehat{X}_{\!/U_\bullet}& \widehat{X}_{\!/W_{\bullet\bullet}}\ar[l]_-{p}\ar[r]^-{q'}&
 \widehat{X}_{\!/V_\bullet}\lefteqn{,}
 }
 \quad
 \xymatrix{%
 \widehat{Z}'_{\!U'_\bullet}\ar[d]^-{f_\alpha}& \widehat{Z}'_{\!W'_{\bullet\bullet}}\ar[l]_-{p'}\ar[d]^-{f_{\alpha\beta}}\ar[r]^-{q'}&
 \widehat{Z}'_{\!V'_\bullet}\ar[d]^-{f_\beta}\\
 \widehat{Z}_{\!/U_\bullet}& \widehat{Z}_{\!/W_{\bullet\bullet}}\ar[l]_-{p}\ar[r]^-{q'}&
 \widehat{Z}_{\!/V_\bullet}\lefteqn{.}
 }
 \]
 Moreover, there is the natural morphism ``$i$'' (resp.\ ``$j$'') from the left (resp.\ right)
 diagram to the special fiber of the middle diagram. These naturally induce the following commutative diagram:
 \[
 \xymatrix@C=15pt{%
  p'^*f_\alpha^*i_U^*Rj_{U*}Rj_U^!R\psi_{\widehat{X}_{\!/U_\bullet}}\Lambda\ar[d]^-{p'^*\widetilde{\varphi}_\alpha}\ar[r]_{\cong}\ar@/^20pt/[rr]^-{f_{\alpha\beta}^*c^*\lambda_{VU}=c'^*f_\red^*\lambda_{VU}}
 &f_{\alpha\beta}^*i_W^*Rj_{W*}Rj_W^!R\psi_{\widehat{X}_{\!/W_{\bullet\bullet}}}\Lambda\ar[d]^-{\widetilde{\varphi}_{\alpha\beta}}
 &q'^*f_\beta^*i_V^*Rj_{V*}Rj_V^!R\psi_{\widehat{X}_{\!/V_\bullet}}\Lambda\ar[d]^-{q'^*\widetilde{\varphi}_\beta}\ar[l]^-{\cong}\\
 p'^*i'^*_{U'}Rj'_{U'*}Rj'^!_{U'}R\psi_{\widehat{X}'_{\!/U'_\bullet}}\Lambda\ar[r]^-{\cong}\ar@/_20pt/[rr]_-{c'^*\lambda'_{V'U'}}
 &i'^*_{W'}Rj'_{W'*}Rj'^!_{W'}R\psi_{\widehat{X}'_{\!/W'_{\bullet\bullet}}}\Lambda
 &q'^*i'^*_{V'}Rj'_{V'*}Rj'^!_V R\psi_{\widehat{X}'_{\!/V'_\bullet}}\Lambda\ar[l]_-{\cong}\lefteqn{.}
 }
 \]
 Therefore the diagram
 \[
  \xymatrix{%
 f_\red^*R\Psi_{\!\mathcal{X},\mathcal{Z},U}\Lambda\ar[rr]^-{f_\red^*\lambda_{VU}}\ar[d]_-{\lambda'^{-1}_{U'}\circ \varphi_\alpha\circ f_\red^*\lambda_U}&& f_\red^*R\Psi_{\!\mathcal{X},\mathcal{Z},V}\Lambda\ar[d]^-{\lambda'^{-1}_{V'}\circ \varphi_\beta\circ f_\red^*\lambda_V}\\
 R\Psi_{\!\mathcal{X}',\mathcal{Z}',U'}\Lambda\ar[rr]^-{\lambda'_{V'U'}}&& R\Psi_{\!\mathcal{X}',\mathcal{Z}',V'}\Lambda
 }
 \]
 is commutative, which implies $\varphi_\alpha=\varphi_\beta$. Thus we can put $\varphi_f=\varphi_\alpha$.
 The compatibility with composition is immediate.
\end{prf}

\begin{rem}
 It is easy to observe that the distinguished triangle in Proposition \ref{prop:excision} is
 compatible with the change of $\mathcal{X}$ in the obvious sense. 
\end{rem}

\begin{prop}\label{prop:open-base-change}
 In the same situation as Proposition \ref{prop:functoriality}, assume that
 $f$ is an open immersion. Then $\varphi_f$ is an isomorphism.
\end{prop}

\begin{prf}
 First we assume that $\mathcal{X}'$ is affine. Let us take an affine open covering $U$ of $\mathcal{X}_\red$.
 It induces the affine open covering $U'$ of $\mathcal{X}'_\red$, since $\mathcal{X}_\red$ is separated.
 Then the canonical morphism $\widehat{X}'_{\!/U'_\bullet}\longrightarrow \widehat{X}_{\!/U_\bullet}$ is regular
 by Proposition \ref{prop:regular}. Therefore the regular base change theorem and \cite[Corollaire 4.7]{Dualite}
 ensure the isomorphy of $\varphi_f$.

 In the general case, let us take an affine open covering $\mathcal{X}'=\bigcup_{i\in I}\mathcal{U}_i$.
 Then for each $i\in I$ we have the following commutative diagram:
 \[
 \xymatrix{%
 (R\Psi_{\!\mathcal{X},\mathcal{Z}}\Lambda)\vert_{(\mathcal{U}_i)_\red}\ar[rr]^-{\varphi_f\vert_{(\mathcal{U}_i)_\red}}\ar[rd]_-{(1)}& &(R\Psi_{\!\mathcal{X}',\mathcal{Z}'}\Lambda)\vert_{(\mathcal{U}_i)_\red}
 \ar[ld]^-{(2)}\\
 & R\Psi_{\mathcal{U}_i,\mathcal{Z}\cap \mathcal{U}_i}\Lambda\lefteqn{.}
 }
 \]
 Since (1) and (2) are isomorphisms, so is $\varphi_f\vert_{(\mathcal{U}_i)_\red}$.
 Therefore $\varphi_f$ is also an isomorphism.
\end{prf}

\begin{rem}
 Later we will prove that $\varphi_f$ is an isomorphism if $f$ is smooth (Proposition \ref{prop:smooth-base-change}).
\end{rem}

\subsubsection{Comparison with nearby cycle for schemes}
\begin{prop}\label{prop:comparison-scheme}
 Let $X$ be a locally noetherian excellent scheme over $S$, $Y$ a closed subscheme of $X_s$ which is separated,
 and $Z=(Z_1,Z_2)$ a pair of closed subschemes of $X_s$ with $Z_1\subset Z_2$. We also denote the locally closed
 subscheme $Z_2\setminus Z_1$ of $X_s$ by $Z$. Let $i\colon Y\hooklongrightarrow X_s$ and 
 $j\colon Z\hooklongrightarrow X_s$ be the natural immersions.
 Let $\mathcal{X}$ be the completion of $X$ along $Y$ and put $\mathcal{Z}_1=\mathcal{X}\times_XZ_1$, 
 $\mathcal{Z}_2=\mathcal{X}\times_XZ_2$ and $\mathcal{Z}=(\mathcal{Z}_1,\mathcal{Z}_2)$.
 Assume $\mathcal{X}$ is excellent.
 Then, there exists a natural isomorphism
 \[
 i^*Rj_*Rj^!R\psi_{X}\Lambda\yrightarrow{\cong} R\Psi_{\!\mathcal{X},\mathcal{Z}}\Lambda,
 \]
 which is functorial on the pair $(X,Y)$.
\end{prop}

\begin{prf}
 We may assume that $Y$ is reduced. Then we can identify $Y$ with $\mathcal{X}_\red$.
 Let us take an affine open covering $V$ of $X$ and $U$ the affine open covering of $Y=\mathcal{X}_\red$
 induced from $V$.
 For an affine open subscheme $\Spec A$ of $X$ belonging to $V$, 
 we have $\widehat{X}_{\!/Y\cap \Spec A}=\Spec \widehat{A}$, where $\widehat{A}$ is the completion of
 $A$ by the defining ideal of $Y\cap \Spec A$. Thus we have the $S$-morphism
 $\widehat{X}_{\!/Y\cap \Spec A}\longrightarrow \Spec A\hooklongrightarrow X$,
 which is regular since $A$ is excellent.
 Therefore, we have the natural augmentation $a\colon \widehat{X}_{\!/U_\bullet}\longrightarrow X$.
 Moreover we have the following commutative diagram whose right rectangle is cartesian:
 \[
 \xymatrix{%
 U_\bullet\ar[r]^-{i_U}\ar[d]^-{a}& (\widehat{X}_{\!/U_\bullet})_s\ar[d]^-{a}& \widehat{Z}_{\!/U_\bullet}\ar[l]_-{j_U}\ar[d]^-{a}\\
 Y\ar[r]^-{i}& X_s& Z\ar[l]_-{j}\lefteqn{.}
 }
 \]
 Therefore we have natural morphisms
 \begin{align*}
  &i^*Rj_*Rj^!R\psi_{X}\Lambda\longrightarrow Ra_*a^*i^*Rj_*Rj^!R\psi_{X}\Lambda
 =Ra_*i_U^*a^*Rj_*Rj^!R\psi_{X}\Lambda\\
  &\qquad \longrightarrow Ra_*i_U^*Rj_{U*}a^*Rj^!R\psi_{X}\Lambda
  \longrightarrow Ra_*i_U^*Rj_{U*}Rj_U^!a^*R\psi_{X}\Lambda\\
  &\qquad \longrightarrow Ra_*i_U^*Rj_{U*}Rj_U^!R\psi_{\widehat{X}_{\!/U_\bullet}}\Lambda
  =R\Psi_{\!\mathcal{X},\mathcal{Z},U}\Lambda
  \yrightarrow[\cong]{\lambda_U}R\Psi_{\!\mathcal{X},\mathcal{Z}}\Lambda,
 \end{align*}
 which are isomorphisms, since $a\colon U_\bullet\longrightarrow Y$ is a morphism of cohomological descent
 and $a_m\colon \widehat{X}_{\!/U_m}\longrightarrow X$ is a regular morphism for every $m\ge 0$.

 By the same method as in the proof of Proposition \ref{prop:functoriality}, we can prove that the morphism
 above is independent of the choice of the affine open covering $V$.
 The functoriality is clear.
\end{prf}

\begin{cor}\label{cor:cor-comp}
 Let the notation be the same as in Proposition \ref{prop:comparison-scheme}.
 Then the object $i^*Rj_*Rj^!R\psi_{X}\Lambda$ of $D^+(Y,\Lambda)$ depends only on 
 $\mathcal{X}$ and $\mathcal{Z}$.
 To give more precise statement, let $X'$, $Y'$ and $Z'=(Z_1',Z_2')$ be another data 
 as in Proposition \ref{prop:comparison-scheme}, and 
 denote the natural immersion $Y'\hooklongrightarrow X'_s$
 (resp.\ $Z'=Z'_2\setminus Z'_1\hooklongrightarrow X'_s$) by $i'$ (resp.\ $j'$).
 Let $\mathcal{X}'$ be the completion of $X'$ along $Y'$ and put 
 $\mathcal{Z}'_1=\mathcal{X}'\times_{X'}Z'_1$, $\mathcal{Z}'_2=\mathcal{X}'\times_{X'}Z'_2$
 and $\mathcal{Z}'=(\mathcal{Z}'_1,\mathcal{Z}'_2)$.
 Assume that the pair $(\mathcal{X}',\mathcal{Z}')$ is isomorphic to $(\mathcal{X},\mathcal{Z})$. 
 Then we have an isomorphism
 $i^*Rj_*Rj^!R\psi_X\Lambda\cong i'^*Rj'_*Rj'^!R\psi_{X'}\Lambda$, under the identification 
 $Y_\red=\mathcal{X}_\red\cong \mathcal{X}'_\red=Y'_\red$.
\end{cor}

\begin{prf}
 Clear from Proposition \ref{prop:comparison-scheme}.
\end{prf}

\begin{rem}
 Corollary \ref{cor:cor-comp} for $Z=(\varnothing,Y)$ is due to Berkovich \cite[Theorem 3.1]{MR1395723}.
\end{rem}

\subsubsection{Smooth base change}
\begin{prop}\label{prop:smooth-base-change}
 In the same situation as Proposition \ref{prop:functoriality}, assume that
 $f$ is smooth. Then $\varphi_f\colon f_\red^*R\Psi_{\!\mathcal{X},\mathcal{Z}}\Lambda\longrightarrow R\Psi_{\!\mathcal{X}',\mathcal{Z}'}\Lambda$ is an isomorphism.
\end{prop}

In order to prove this proposition, we need the following lemma:

\begin{lem}\label{lem:etale-algebraization}
 Let $A$ be a noetherian ring and $I$ an ideal of $A$. Assume that $A$ is $I$-adically complete.
 Let $f_1,\ldots,f_m\in A\langle T_1,\ldots,T_m\rangle$, where $A\langle T_1,\ldots,T_m\rangle$ is
 the $I$-adic completion of $A[T_1,\ldots,T_m]$, and put $B=A\langle T_1,\ldots,T_m\rangle/(f_1\ldots,f_m)$.
 If the image of $\Delta=\det (\partial f_i/\partial T_j)_{i,j}$ in $B$ is invertible,
 then there exists an \'etale $A$-algebra whose $I$-adic completion is isomorphic to $B$ as an $A$-algebra.
\end{lem}

\begin{prf}
 We can take $g_1,\ldots,g_m\in A[T_1,\ldots,T_m]$ such that $g_i-f_i\in IA\langle T_1,\ldots,T_m\rangle$.
 Put $B'=A\langle T_1,\ldots,T_m\rangle/(g_1,\ldots,g_m)$. We will prove
 $B\cong B'$ as $A$-algebras.

 Put $\Delta'=\det (\partial g_i/\partial T_j)_{i,j}$. Clearly $\Delta\equiv \Delta'\pmod{I}$.
 Therefore the image of $\Delta$ and $\Delta'$ in $B/IB$ are the same.
 On the other hand, we have $B/IB\cong B'/IB'$. Thus the image of $\Delta'$ in $B'/IB'$ is invertible.
 Since $B'$ is $IB'$-adically complete, this implies that the image of $\Delta'$ in $B'$ is invertible.
 Therefore $\Spf B$ and $\Spf B'$ are \'etale over $\Spf A$ and their restriction to $\Spf A/I$
 are isomorphic. Hence we have $\Spf B\cong \Spf B'$ and $B\cong B'$.

 Put $C=A[T_1,\ldots,T_m,T_{m+1}]/(g_1,\ldots,g_m,T_{m+1}\Delta'-1)$. Clearly $C$ is \'etale over $A$.
 On the other hand, since $\Delta'$ is invertible in $B'/IB'$, 
 the $I$-adic completion of $C$ is isomorphic to $B'$ as an $A$-algebra.
 This completes the proof.
\end{prf}

\begin{prf}[Proposition \ref{prop:smooth-base-change}]
 By Proposition \ref{prop:open-base-change}, we may assume that $\mathcal{X}$ is an affine formal scheme $\Spf A$.
 Moreover, we may assume that there exists an integer $m\ge 0$ and an \'etale morphism
 $g\colon \mathcal{X}'\longrightarrow \Spf A\langle T_1,\ldots,T_m\rangle$ such that the composite
 $\mathcal{X}'\yrightarrow{g}\Spf A\langle T_1,\ldots,T_m\rangle\yrightarrow{\pr_m} \Spf A=\mathcal{X}$ is
 equal to $f$.
 Thus we have only to consider the following two cases:
 \begin{itemize}
  \item $f=\pr_m$.
  \item $f$ is \'etale.
 \end{itemize}
 First consider the case where $f=\pr_m$.
 Let us denote the defining ideal of $\mathcal{X}_\red$, $\mathcal{Z}_1$, $\mathcal{Z}_2$ by $I$, $J_1$, $J_2$
 respectively, and put $X=\Spec A$, $Y=\Spec A/I$, $Z_1=\Spec A/J_1$, $Z_2=\Spec A/J_2$, $Z=Z_2\setminus Z_1$,
 $X'=\Spec A[T_1,\ldots,T_m]$, $Y'=X'\times_XY$, $Z'_1=X'\times_XZ_1$, $Z'_2=X'\times_XZ_2$
 and $Z'=Z'_2\setminus Z'_1$.
 Denote the natural immersions $Y\hooklongrightarrow X_s$, $Z\hooklongrightarrow X_s$,
 $Y'\hooklongrightarrow X'_s$, $Z'\hooklongrightarrow X'_s$ by $i$, $j$, $i'$, $j'$, respectively.
 
 Then, by Proposition \ref{prop:comparison-scheme} (or the definition of $R\Psi_{\!\mathcal{X},\mathcal{Z}}\Lambda$), we have the isomorphisms
 $i^*Rj_*Rj^!R\psi_X\Lambda\yrightarrow{\cong}R\Psi_{\!\mathcal{X},\mathcal{Z}}\Lambda$ and
 $i'^*Rj'_*Rj'^!R\psi_{X'}\Lambda\yrightarrow{\cong}R\Psi_{\!\mathcal{X}',\mathcal{Z}'}\Lambda$.
 Furthermore, by the functoriality of the isomorphism in Proposition \ref{prop:comparison-scheme},
 we have the following commutative diagram:
 \[
  \xymatrix{%
 (\pr_m)_\red^*i^*Rj_*Rj^!R\psi_X\Lambda\ar[r]^-{\cong}\ar[d]& 
 (\pr_m)_\red^*R\Psi_{\!\mathcal{X},\mathcal{Z}}\Lambda\ar[d]^-{\varphi_f}\\
 i'^*Rj'_*Rj'^!R\psi_{X'}\Lambda\ar[r]^-{\cong}& R\Psi_{\!\mathcal{X}',\mathcal{Z}'}\Lambda\lefteqn{.}
 }
 \]
 Thus it suffices to show that the natural morphism 
 $(\pr_m)_\red^*i^*Rj_*Rj^!R\psi_X\Lambda\longrightarrow i'^*Rj'_*Rj'^!R\psi_{X'}\Lambda$ is an isomorphism.
 This is an immediate corollary of the smooth base change theorem.

 Next we consider the case where $f$ is \'etale. By shrinking $\mathcal{X}'$ if necessary, 
 we may assume that $\mathcal{X}'=\Spf B$, where $B=A\langle T_1,\ldots,T_m\rangle/(f_1,\ldots,f_m)$
 such that the image of $\det (\partial f_i/\partial T_j)_{i,j}$ in $B$ is invertible.
 By Lemma \ref{lem:etale-algebraization}, there exists an \'etale $A$-algebra $C$ such that
 the $I$-adic completion of $C$ is isomorphic to $B'$ as an $A$-algebra.
 Let $X$, $Y$, $Z$ be the same as above. Put $X'=\Spec C$ and define $Y'$ and $Z'$ as above.
 Then, by exactly the same way as above, we can prove that $\varphi_f$ is an isomorphism.
\end{prf}

\subsubsection{Invariance under admissible blow-up}
\begin{prop}\label{prop:inv-adm-bu}
 Let $\pi\colon \mathcal{X}'\longrightarrow \mathcal{X}$ be
 an admissible blow-up \cite[Definition 4.1.1]{MR1360610}. 
 Let $\mathcal{Z}=(\mathcal{Z}_1,\mathcal{Z}_2)$ be a pair of closed formal subschemes of $\mathcal{X}$
 with $\mathcal{Z}_1\subset \mathcal{Z}_2$ and put $\mathcal{Z}'=(\mathcal{X}'\times_{\mathcal{X}}\mathcal{Z}_1,\mathcal{X}'\times_{\mathcal{X}}\mathcal{Z}_2)$.
 Assume moreover that $\mathcal{X}'$ is excellent.
 Then the natural morphism $R\Psi_{\!\mathcal{X},\mathcal{Z}}\Lambda\longrightarrow R\pi_*R\Psi_{\!\mathcal{X}',\mathcal{Z}'}\Lambda$ induced by $\varphi_\pi$ is an isomorphism.
\end{prop}

\begin{prf}
 By Proposition \ref{prop:open-base-change}, we may assume that $\mathcal{X}$ is affine.
 Let $\mathcal{X}=\Spf A$, $\mathcal{X}_\red=\Spf A/I$, $\mathcal{Z}_1=\Spf A/J_1$ and $\mathcal{Z}_2=\Spf A/J_2$.
 Put $X=\Spec A$, $Y=\Spec A/I$, $Z_1=\Spec A/J_1$, $Z_2=\Spec A/J_2$, $Z=Z_2\setminus Z_1$ and 
 let $i\colon Y\hooklongrightarrow X_s$ and $j\colon Z\hooklongrightarrow X_s$ be the natural immersions.
 Then, by Proposition \ref{prop:comparison-scheme} (or the definition of $R\Psi_{\!\mathcal{X},\mathcal{Z}}\Lambda$), we have the isomorphism
 $i^*Rj_*Rj^!R\psi_X\Lambda\yrightarrow{\cong}R\Psi_{\!\mathcal{X},\mathcal{Z}}\Lambda$.

 On the other hand, by the definition of an admissible blow-up, there exists a blow-up
 $\widetilde{\pi}\colon X'\longrightarrow X$ whose center is contained in $Y$ such that
 the $I$-adic completion of it is isomorphic to $\pi\colon \mathcal{X}'\longrightarrow \mathcal{X}$.
 Therefore, if we put $Y'=X'\times_XY$, $Z'_1=X'\times_XZ_1$, $Z'_2=X'\times_XZ_2$, $Z'=X'\times_XZ$
 and denote the natural immersion $Y'\hooklongrightarrow X'_s$ (resp.\ $Z'\hooklongrightarrow X'_s$)
 by $i'$ (resp.\ $j'$), we have the isomorphism
 $i'^*Rj'_*Rj'^!R\psi_{X'}\Lambda\yrightarrow{\cong}R\Psi_{\!\mathcal{X}',\mathcal{Z}'}\Lambda$.
 Furthermore, by the functoriality of the isomorphism in Proposition \ref{prop:comparison-scheme},
 we have the following commutative diagram:
 \[
 \xymatrix{%
 i^*Rj_*Rj^!R\psi_X\Lambda\ar[r]^-{\cong}\ar[d]& 
 R\Psi_{\!\mathcal{X},\mathcal{Z}}\Lambda\ar[d]\\
 \pi_{\red *}i'^*Rj'_*Rj'^!R\psi_{X'}\Lambda\ar[r]^-{\cong}& \pi_{\red *}R\Psi_{\!\mathcal{X}',\mathcal{Z}'}\Lambda\lefteqn{.}
 }
 \]
 Therefore it suffices to show that the natural morphism
 \[
  i^*Rj_*Rj^!R\psi_X\Lambda\yrightarrow{(*)} \pi_{\red *}i'^*Rj'_*Rj'^!R\psi_{X'}\Lambda
 \]
 is an isomorphism. By the proper base change theorem, we have
 \begin{align*}
  \pi_{\red *}i'^*Rj'_*Rj'^!R\psi_{X'}\Lambda&\yleftarrow{\cong}i^*\widetilde{\pi}_*Rj'_*Rj'^!R\psi_{X'}\Lambda
 =i^*Rj_*\widetilde{\pi}_*Rj'^!R\psi_{X'}\Lambda\\
  &\cong i^*Rj_*Rj^!\widetilde{\pi}_*R\psi_{X'}\Lambda\yleftarrow{\cong}i^*Rj_*Rj^!R\psi_X\widetilde{\pi}_*\Lambda=i^*Rj_*Rj^!R\psi_X\Lambda.
 \end{align*}
 It is not difficult to show that the isomorphism above coincides with the morphism $(*)$.
 This completes the proof.
\end{prf}

\begin{rem}
 It is naturally expected that the excellence of $\mathcal{X}'$ automatically
 follows from the excellence of $\mathcal{X}$.
\end{rem}

\subsubsection{Finiteness}
\begin{defn}\label{defn:loc-alg}
 Let $\mathcal{X}$ be a locally noetherian formal scheme over $\mathcal{S}$.
 We say that $\mathcal{X}$ is algebraizable
 if there exist a scheme $X$ which is separated and locally of finite type over $S$ and
 a closed subscheme $Y$ of $X_s$ such that
 the completion of $X$ along $Y$ is isomorphic to $\mathcal{X}$ over $\mathcal{S}$.
 We say that $\mathcal{X}$ is locally algebraizable
 if there exists an open covering $\mathcal{X}=\bigcup_{i\in I}\mathcal{U}_i$ such that
 $\mathcal{U}_i$ is algebraizable for each $i\in I$.
\end{defn}

\begin{rem}
 By Proposition \ref{prop:special-exc}, a locally algebraizable formal scheme is automatically excellent.
\end{rem}

\begin{prop}\label{prop:finiteness}
 Assume that $\mathcal{X}$ is locally algebraizable. 
 Then $R\Psi_{\!\mathcal{X},\mathcal{Z}}\Lambda$ is a constructible complex
 for every pair $\mathcal{Z}=(\mathcal{Z}_1,\mathcal{Z}_2)$ of closed formal subschemes of $\mathcal{X}$
 with $\mathcal{Z}_1\subset \mathcal{Z}_2$.
\end{prop}

\begin{prf}
 By Proposition \ref{prop:open-base-change}, we may assume that 
 there exist an affine scheme $X=\Spec A$ of finite type over $S$ and a reduced closed subscheme
 $Y=\Spec A/I$ of $X_s$ such that $\mathcal{X}$ is the completion of $X$ along $Y$.
 Let $\widehat{A}$ be the $I$-adic completion of $A$. Then $\mathcal{X}=\Spf \widehat{A}$. 
 Let $\mathcal{Z}_1=\Spf \widehat{A}/J_1$, $\mathcal{Z}_2=\Spf \widehat{A}/J_2$ and
 put $\widehat{Z}_1=\Spec \widehat{A}/J_1$, $\widehat{Z}_2=\Spec \widehat{A}/J_2$ and
 $\widehat{Z}=\widehat{Z}_2\setminus \widehat{Z}_1$.
 By the construction, 
 $R\Psi_{\!\mathcal{X},\mathcal{Z}}\Lambda\cong i^*Rj_*Rj^!R\psi_{\widehat{X}}\Lambda$,
 where $i\colon Y\hooklongrightarrow \widehat{X}_s$ and $j\colon \widehat{Z}\hooklongrightarrow \widehat{X}_s$
 denote the natural immersions.

 Let us denote by $h$ the natural morphism $\widehat{X}\longrightarrow X$, which is regular since $A$
 is excellent. Then, by the regular base change theorem, we have 
 $R\psi_{\widehat{X}}\Lambda\cong h^*R\psi_X\Lambda$. Therefore $R\psi_{\widehat{X}}\Lambda$ is constructible,
 since $R\psi_X\Lambda$ is constructible by \cite[Finitude, Th\'eor\`eme 3.2]{MR0463174}.
 Let us factorize $j$ as $\widehat{Z}\stackrel{j_1}{\hooklongrightarrow}\widehat{Z}_2\stackrel{j_2}{\hooklongrightarrow}\widehat{X}_s$, where $j_1$ is an open immersion and $j_2$ is a closed immersion.
 Denote the open immersion $\widehat{X}_s\setminus \widehat{Z}_2\hooklongrightarrow \widehat{X}_s$
 by $j'$.
 Consider the distinguished triangle
 \[
  Rj_2^!R\psi_{\widehat{X}}\Lambda\longrightarrow j_2^*R\psi_{\widehat{X}}\Lambda\longrightarrow
 j_2^*Rj'_*j'^*R\psi_{\widehat{X}}\Lambda\longrightarrow Rj_2^!R\psi_{\widehat{X}}\Lambda[1].
 \]
 We know that $j_2^*R\psi_{\widehat{X}}\Lambda$ is constructible.
 Moreover, by Gabber's finiteness theorem \cite{Gabber-finiteness},
 $j_2^*Rj'_*j'^*R\psi_{\widehat{X}}\Lambda$ is also constructible
 (note that $\widehat{X}_s$ is excellent and $j'$ is quasi-compact).
 Therefore $Rj_2^!R\psi_{\widehat{X}}\Lambda$ is constructible, and thus 
 $Rj^!R\psi_{\widehat{X}}\Lambda=Rj_1^*Rj_2^!R\psi_{\widehat{X}}\Lambda$ is constructible.
 Finally, by Gabber's theorem again, $i^*Rj_*Rj^!R\psi_{\widehat{X}}\Lambda$ is also constructible.
 This completes the proof.
\end{prf}

\begin{rem}
 By using the $\ell$-adic formalism in \cite{MR1106899}, We can also define 
 $R\Psi_{\!\mathcal{X},\mathcal{Z}}\Lambda$ for $\Lambda=\Z_\ell$ or $\Q_\ell$.
 All the properties in this subsection hold for these $\ell$-adic coefficients.
\end{rem}

\section{Comparison with rigid geometry}
\subsection{Comparison with Berkovich's formal nearby cycle}
In this subsection, we will compare $R\Psi_{\!\mathcal{X}}\Lambda$ with the formal nearby cycle of Berkovich.
First we recall the definition of the formal nearby cycle functor.
We will use the framework of adic spaces due to Huber \cite{MR1306024}.
By the dictionary in \cite[8.3]{MR1734903},
we can see without difficulty that our definition coincides with that in \cite{MR1395723}.

Let $\mathcal{X}$ be a locally noetherian formal scheme over $\mathcal{S}$. 

\begin{defn}
 \begin{enumerate}
  \item We put $t(\mathcal{X})_\eta=t(\mathcal{X})\setminus V(\varpi)$;
	recall that $\varpi$ is a uniformizer of $R$.
	It is an adic space over $\Spa(F,R)$, where $F=\Frac R$.
  \item We denote the category of admissible blow-ups of $\mathcal{X}$ by $\Phi_\mathcal{X}$.
	For an object $\mathcal{X}'\longrightarrow \mathcal{X}$, $\mathcal{X}'_\varpi$ denotes 
	the maximal open formal subscheme of $\mathcal{X}'$ where $\varpi\mathcal{O}_{\mathcal{X}'}$ is
	an ideal of definition of $\mathcal{X}'$. Obviously $\mathcal{X}'_\varpi$ is a $\varpi$-adic formal scheme.
 \end{enumerate}
\end{defn}

\begin{rem}
 The adic space $t(\mathcal{X})_\eta$ is denoted by $\widetilde{d}(\mathcal{X})$ in \cite[3.12]{MR1620118}.
 If $\mathcal{X}$ is $\varpi$-adic, then it coincides with $d(\mathcal{X})$ in \cite[1.9]{MR1734903},
 which is quasi-compact if $\mathcal{X}$ is quasi-compact.
\end{rem}

\begin{lem}\label{lem:generic-fiber}
 \begin{enumerate}
  \item Let $\mathcal{X}'\longrightarrow \mathcal{X}$ be an object of $\Phi_{\mathcal{X}}$.
	Then the natural morphism $t(\mathcal{X}')_\eta\longrightarrow t(\mathcal{X})_\eta$ is an isomorphism.
  \item Let $(\mathcal{X}''\longrightarrow \mathcal{X})\longrightarrow (\mathcal{X}'\longrightarrow \mathcal{X})$
	be a morphism in $\Phi_{\mathcal{X}}$. Then an open immersion
	$t(\mathcal{X}'_\varpi)_\eta\hooklongrightarrow t(\mathcal{X}''_\varpi)_\eta$ is naturally induced.
  \item We have $t(\mathcal{X})_\eta=\varinjlim_{(\mathcal{X}'\to \mathcal{X})\in \Phi^{\mathrm{op}}_{\mathcal{X}}}t(\mathcal{X}'_\varpi)_\eta$.
  \item If $\mathcal{X}$ is special, then $t(\mathcal{X})_\eta$ is locally of finite type over $\Spa (F,R)$.
  \item For special formal schemes $\mathcal{X}$ and $\mathcal{Y}$ over $\mathcal{S}$,
	$\mathcal{X}\times_{\mathcal{S}}\mathcal{Y}$ is also special and
	we have $t(\mathcal{X}\times_{\mathcal{S}}\mathcal{Y})_\eta=t(\mathcal{X})_\eta\times_{\Spa(F,R)}t(\mathcal{Y})_\eta$.
 \end{enumerate}
\end{lem}

\begin{prf}
 i) is well-known. For ii), let $\mathcal{U}''$ be the inverse image of $\mathcal{X}'_\varpi$ in 
 $\mathcal{X}''$. Clearly $\mathcal{U}''\subset \mathcal{X}''_\varpi$. On the other hand,
 the natural morphism $\mathcal{U}''\longrightarrow \mathcal{X}'_\varpi$ is an admissible blow-up.
 Thus by i) we have the open immersion $t(\mathcal{X}'_\varpi)_\eta\yleftarrow{\cong}t(\mathcal{U}'')_\eta\hooklongrightarrow t(\mathcal{X}''_\varpi)_\eta$. It is easy to observe that this gives the contravariant functor
 from $\Phi_{\mathcal{X}}$ to the category $\mathbf{Qc}_{t(\mathcal{X})_\eta}$
 of quasi-compact open adic subspaces of $t(\mathcal{X})_\eta$.

 We prove iii). 
 Since $t(\mathcal{X})_\eta=\varinjlim_{\mathbb{U}\in \mathbf{Qc}_{t(\mathcal{X})_\eta}}\mathbb{U}$,
 it suffices to show that the functor 
 $\Phi_{\mathcal{X}}^{\mathrm{op}}\longrightarrow \mathbf{Qc}_{t(\mathcal{X})_\eta}$ introduced above is 
 cofinal. Since $\Phi^{\mathrm{op}}_{\mathcal{X}}$ is filtered and $\mathbf{Qc}_{t(\mathcal{X})_\eta}$ is a filtered ordered
 set, it suffices to show that for every $\mathbb{U}\in \mathbf{Qc}_{t(\mathcal{X})_\eta}$
 there exists an object $\mathcal{X}'\longrightarrow \mathcal{X}$ of $\Phi_{\mathcal{X}}$ such that
 $\mathbb{U}\subset t(\mathcal{X}'_\varpi)_\eta$.
 We may assume that $\mathcal{X}$ is an affine formal scheme $\Spf A$ (use \cite[I, (9.4.7)]{EGA} to
 extend a local admissible blow-up to global). Let $I=(a_1,\ldots,a_n)$ be an ideal of definition of $A$.
 Then it is easy to see that
 \[
  t(\Spf A)\setminus V(\varpi)=\bigcup_{k\ge 1}R\Bigl(\dfrac{a_1^k,\ldots,a_n^k}{\varpi}\Bigr),
 \]
 where
 \[
 R\Bigl(\dfrac{a_1^k,\ldots,a_n^k}{\varpi}\Bigr)
 =\bigl\{v\in \Spa (A,A)\bigm| v(a_i^k)\le v(\varpi)\neq 0\ (i=1,\ldots,n)\bigr\}
 \]
 is a rational subset of $\Spa (A,A)$. Let $\mathcal{X}'\longrightarrow \mathcal{X}=\Spa A$ be the
 admissible blow-up along the open ideal $(\varpi,a_1^k,\ldots,a_n^k)$ of $A$. 
 Then we have $t(\mathcal{X}'_\varpi)_\eta=R\Bigl(\dfrac{a_1^k,\ldots,a_n^k}{\varpi}\Bigr)$.
 This completes the proof of iii).

 For iv), note that $t(\mathcal{X}'_\varpi)_\eta$ is an adic space of finite type over $\Spa (F,R)$
 for every object $\mathcal{X}'\longrightarrow \mathcal{X}$ of $\Phi_{\mathcal{X}}$.

 We prove v). By \cite[Lemma 1.1 (v)]{MR1395723}, $\mathcal{X}\times_{\mathcal{S}}\mathcal{Y}$ is also special.
 To observe $t(\mathcal{X}\times_\mathcal{S}\mathcal{Y})_\eta=t(\mathcal{X})_\eta\times_{\Spa(F,R)}t(\mathcal{Y})_\eta$, it suffices to show the equality
 $t(\mathcal{X}\times_\mathcal{S}\mathcal{Y})=t(\mathcal{X})\times_{t(\mathcal{S})}t(\mathcal{Y})$.
 We may assume that $\mathcal{X}$ and
 $\mathcal{Y}$ are affine. By \cite[(1) in the proof of Proposition 4.1]{MR1306024},
 for an adic space $\mathbb{X}$ and a commutative diagram of locally and topologically ringed spaces
 \[
  \xymatrix{%
 (\mathbb{X},\mathcal{O}^+_{\mathbb{X}})\ar[r]\ar[d]& \mathcal{X}\ar[d]\\
 \mathcal{Y}\ar[r]& \mathcal{S}\lefteqn{,}
 }
 \]
 we have a unique morphism $(\mathbb{X},\mathcal{O}^+_{\mathbb{X}})\longrightarrow \mathcal{X}\times_{\mathcal{S}}\mathcal{Y}$ that makes the obvious diagram commutative. 
 It implies $t(\mathcal{X}\times_{\mathcal{S}}\mathcal{Y})=t(\mathcal{X})\times_{t(\mathcal{S})}t(\mathcal{Y})$,
 which concludes the proof.
\end{prf}

\begin{defn}
 Let $\overline{R}$ be the ring of integers of $\overline{F}$, the fixed separable closure of $F$.
 For a special formal scheme $\mathcal{X}$ over $\mathcal{S}$, we put
 \[
  t(\mathcal{X})_{\overline{\eta}}=t(\mathcal{X})_\eta\times_{\Spa(F,R)}\Spa(\overline{F},\overline{R}).
 \]
 By Lemma \ref{lem:generic-fiber} iv) and \cite[Proposition 3.7]{MR1306024}, the fiber product above
 can be defined.
\end{defn}

Now we can define the formal nearby cycle functor. 

\begin{defn}
 Let $\mathcal{X}$ be a special formal scheme over $\mathcal{S}$.
 Consider the functor $(\mathcal{X}_\red)_\et\cong \mathcal{X}_\et\longrightarrow (t(\mathcal{X})_{\overline{\eta}})_\et$; $\mathcal{Y}\longmapsto t(\mathcal{Y})_{\overline{\eta}}$. By \cite[Lemma 3.5.1]{MR1734903} and
 \cite[Lemma 3.9 (i)]{MR1306024}, this gives a morphism of sites
 $\Psi^\ad_{\!\mathcal{X}}\colon (t(\mathcal{X})_{\overline{\eta}})_\et\longrightarrow (\mathcal{X}_\red)_\et$.
 We denote the derived functor of $(\Psi^\ad_{\!\mathcal{X}})_*$ by $R\Psi^\ad_{\!\mathcal{X}}$ and 
 call it the formal nearby cycle functor.
\end{defn}

The following is our main comparison result in this subsection.

\begin{thm}\label{thm:comp-Berkovich}
 For a special formal scheme $\mathcal{X}$ over $\mathcal{S}$ such that $\mathcal{X}_\red$ is separated,
 we have a natural morphism
 $\varepsilon^*\colon R\Psi_{\!\mathcal{X}}\Lambda\longrightarrow R\Psi^\ad_{\!\mathcal{X}}\Lambda$.
 Moreover, if $\mathcal{X}$ is locally algebraizable or $\varpi$-adic, then $\varepsilon^*$ is an isomorphism.
\end{thm}

First consider the case where $\mathcal{X}$ is an affine formal scheme $\Spf A$.
Let $I$ be the ideal of definition of $A$ such that $A/I$ is reduced, and put $\widehat{X}=\Spec A$, $Y=\Spec A/I$.
We will define the morphism of sites $(t(\mathcal{X})_{\overline{\eta}})_\et\longrightarrow (\widehat{X}_{\overline{\eta}})_\et$. 

\begin{lem}
 Let $Z$ be an adic space locally of finite type over $\Spa (F,R)$, $T$ a scheme over $S$
 and $(Z,\mathcal{O}_Z)\longrightarrow T$ a morphism of locally ringed spaces over $S$.
 Then we have a natural morphism of locally ringed spaces
 $(Z_{\overline{\eta}},\mathcal{O}_{Z_{\overline{\eta}}})\longrightarrow T_{\overline{\eta}}$.
  
\end{lem}

\begin{prf}
 We may assume that $T=\Spec B$ is affine.
 Take an affinoid open subspace $\Spa(C,C^+)$ of $Z$ where $(C,C^+)$ is complete,
 and construct the morphism 
 \[
  \Spa (C,C^+)\times_{\Spa (F,R)}\Spa (\overline{F},\overline{R})\longrightarrow \Spec (B\otimes_R\overline{F}).
 \]
 The given morphism $Z\longrightarrow T$ induces a morphism $\Spa (C,C^+)\longrightarrow \Spec B$.
 This corresponds to an $R$-homomorphism $B\longrightarrow C$. Since $\varpi$ is invertible in $C$, 
 we get the $F$-homomorphism $B\otimes_RF\longrightarrow C$ and 
 the $\overline{F}$-homomorphism $B\otimes_R\overline{F}\longrightarrow C\otimes_F\overline{F}$.

 On the other hand, by the construction of the fiber product \cite[Proposition 3.7]{MR1306024}, 
 $\Spa (C,C^+)\times_{\Spa (F,R)}\Spa (\overline{F},\overline{R})$ is an affinoid $\Spa (D,D^+)$,
 where $(D,D^+)$ is complete. 
 In particular, we have the following diagram of rings:
 \[
  \xymatrix{%
 D& C\ar[l]\\ \overline{F}\ar[u]& F\lefteqn{.}\ar[l]\ar[u]
 }
 \]
 Thus we have a natural $\overline{F}$-homomorphism $C\otimes_F\overline{F}\longrightarrow D$,
 which corresponds to $\Spa (C,C^+)\times_{\Spa (F,R)}\Spa (\overline{F},\overline{R})=\Spa (D,D^+)\longrightarrow \Spec (C\otimes_F\overline{F})$.
 By composing it with $\Spec (C\otimes_F\overline{F})\longrightarrow \Spec (B\otimes_R\overline{F})$,
 we get the desired morphism of locally ringed spaces.
\end{prf}

By applying this lemma to the natural morphism 
$t(\mathcal{X})_{\eta}\longrightarrow t(\mathcal{X})\longrightarrow \widehat{X}$,
we have the morphism of locally ringed spaces
$t(\mathcal{X})_{\overline{\eta}}\longrightarrow \widehat{X}_{\overline{\eta}}$.
Therefore, for a morphism $W\longrightarrow \widehat{X}_{\overline{\eta}}$ locally of finite type,
we can form the fiber product $t(\mathcal{X})_{\overline{\eta}}\times_{\widehat{X}_{\overline{\eta}}}W$
in the sense of \cite[Proposition 3.8]{MR1306024}.
The functor $W\longmapsto t(\mathcal{X})_{\overline{\eta}}\times_{\widehat{X}_{\overline{\eta}}}W$ gives
a morphism of sites $\varepsilon\colon (t(\mathcal{X})_{\overline{\eta}})_\et\longrightarrow (\widehat{X}_{\overline{\eta}})_\et$ \cite[3.2.8]{MR1734903}.

\begin{lem}\label{lem:natural-morphism}
 Consider the following diagram of sites, where $i$ and $\overline{\jmath}$ denote the natural morphisms:
 \[
  \xymatrix{%
 \bigl(t(\mathcal{X})_{\overline{\eta}}\bigr)_\et\ar[r]^-{\Psi^\ad_{\!\mathcal{X}}}\ar[d]^-{\varepsilon}& Y_\et\ar[d]^-{i}\\
 (\widehat{X}_{\overline{\eta}})_\et\ar[r]^-{\overline{\jmath}}& \widehat{X}_\et\lefteqn{.}
 }
 \]
 We have a natural morphism of functors $(\Psi^\ad_{\!\mathcal{X}})^{-1}\circ i^{-1}\longrightarrow \varepsilon^{-1}\circ \overline{\jmath}^{-1}$.
\end{lem}

\begin{prf}
 Let $W\longrightarrow \widehat{X}$ be an \'etale morphism and denote the $I$-adic completion of $W$ by $\mathcal{W}$.
 Then $(\Psi^\ad_{\!\mathcal{X}})^{-1}(i^{-1}(W))=t(\mathcal{W})_{\overline{\eta}}$.
 On the other hand, 
 $\varepsilon^{-1}(\overline{\jmath}^{-1}(W))=t(\mathcal{X})_{\overline{\eta}}\times_{\widehat{X}_{\overline{\eta}}}W_{\overline{\eta}}$.
 Since we have a natural morphism of locally ringed spaces $t(\mathcal{W})_{\overline{\eta}}\longrightarrow W_{\overline{\eta}}$ by Lemma \ref{lem:natural-morphism}, we have the morphism
 $t(\mathcal{W})_{\overline{\eta}}\longrightarrow t(\mathcal{X})_{\overline{\eta}}\times_{\widehat{X}_{\overline{\eta}}}W_{\overline{\eta}}$ of adic spaces.
\end{prf}

\begin{defn}
 By the lemma above, we have a natural morphism
 \[
 R\Psi_{\!\mathcal{X}}\Lambda=(R\psi_{\widehat{X}}\Lambda)\vert_Y=i^*R\overline{\jmath}_*\Lambda\longrightarrow R(\Psi^\ad_{\!\mathcal{X}})_*\varepsilon^*\Lambda
 =R\Psi^\ad_{\!\mathcal{X}}\Lambda.
 \]
 We denote it by $\varepsilon^*$.
\end{defn}

\begin{prop}\label{prop:affine-alg}
 Assume moreover that there exist a finitely generated $R$-algebra $A_0$
 and an ideal $I_0$ of $A_0$ such that the $I_0$-adic completion of $A_0$ is isomorphic to $A$ as
 a topological $R$-algebra. Then $\varepsilon^*$ is an isomorphism.
\end{prop}

\begin{prf}
 Put $X=\Spec A_0$. By replacing $I_0$ by $\sqrt{I_0}$, we may assume that $\Spec A_0/I_0=Y$.
 Then we have the following commutative diagram:
 \[
  \xymatrix{%
 & (R\psi_X\Lambda)\vert_Y\ar[ld]_-{(1)}^-{\cong}\ar[rd]^-{(2)}_-{\cong}\\
 R\Psi_{\!\mathcal{X}}\Lambda\ar[rr]^-{\varepsilon^*}&&R\Psi^\ad_{\!\mathcal{X}}\Lambda\lefteqn{.}
 }
 \]
 Here (1) is the morphism in Proposition \ref{prop:comparison-scheme}, which is an isomorphism.
 The morphism (2) can be constructed in the similar way as $\varepsilon^*$.
 It is proved to be an isomorphism by \cite[Theorem 3.1]{MR1395723}. Thus $\varepsilon^*$ is also
 an isomorphism.
\end{prf}

\begin{prop}\label{prop:affine-pi-adic}
 Assume that the topology of $A$ is $\varpi$-adic. Then $\varepsilon^*$ is an isomorphism.
\end{prop}

\begin{prf}
 This is a special case of \cite[Theorem 3.5.13]{MR1734903}; apply it to 
 \[
  Y'\stackrel{i'}{\hooklongrightarrow}\widehat{X}\otimes_{R}\overline{R}\stackrel{j'}{\hooklongleftarrow}
 (\widehat{X}\otimes_{R}\overline{R})\setminus Y'=\widehat{X}_{\overline{\eta}},
 \]
 where $Y'$ is the closed subscheme of 
 $\widehat{X}\otimes_{R}\overline{R}$ defined by $\varpi$.
 Note that $R\psi_{\widehat{X}}\Lambda=i'^*Rj'_*\Lambda$,
 since $Y'\longrightarrow Y$ induces an isomorphism on their \'etale sites.
\end{prf}

Now the proof of Theorem \ref{thm:comp-Berkovich} is easy:

\begin{prf}[Theorem \ref{thm:comp-Berkovich}]
 First we will construct the morphism 
 $\varepsilon^*\colon R\Psi_{\!\mathcal{X}}\Lambda\longrightarrow R\Psi^\ad_{\!\mathcal{X}}\Lambda$.
 Take an affine open covering $U=\{U_i\}$ of $\mathcal{X}_\red$ and denote the associated hypercovering
 by $a\colon U_\bullet\longrightarrow \mathcal{X}_\red$.
 Then, for every $i$, $t(\mathcal{X}_{/U_i})_{\overline{\eta}}$ is an open adic subspace of 
 $t(\mathcal{X})_{\overline{\eta}}$ and $\{t(\mathcal{X}_{/U_i})_{\overline{\eta}}\}$ forms an open covering
 of $t(\mathcal{X})_{\overline{\eta}}$. We denote the hypercovering associated with this open covering
 by $a^\ad\colon t(\mathcal{X}_{/U_\bullet})_{\overline{\eta}}\longrightarrow t(\mathcal{X})_{\overline{\eta}}$.
 Then we have the following diagram of sites:
 \[
  \xymatrix{%
 \bigl(t(\mathcal{X})_{\overline{\eta}}\bigr)_\et\ar[r]^-{\Psi^\ad_{\!\mathcal{X}}}& (\mathcal{X}_\red)_\et\\
 \bigl(t(\mathcal{X}_{/U_\bullet})_{\overline{\eta}}\bigr)_\et\ar[u]_-{a^\ad}\ar[d]^-{\varepsilon}\ar[r]^-{\Psi^\ad_{\!\mathcal{X}_{/U_\bullet}}}& (U_\bullet)_\et\ar[d]^-{i}\ar[u]_-{a}\\
 \bigl((\widehat{X}_{\!/U_\bullet})_{\overline{\eta}}\bigr)_\et\ar[r]^-{\overline{\jmath}}& (\widehat{X}_{\!/U_\bullet})_\et\lefteqn{.}
 }
 \]
 Here the lower rectangle is constructed in the same way as the diagram in Lemma \ref{lem:natural-morphism}.
 It is not commutative, but we have a natural morphism of functors
 $(\Psi^\ad_{\!\mathcal{X}_{/U_\bullet}})^{-1}\circ i^{-1}\longrightarrow \varepsilon^{-1}\circ \overline{\jmath}^{-1}$.
 On the other hand, the upper rectangle is obviously commutative.

 Thus, we can construct the morphism $\varepsilon^*\colon R\Psi_{\!\mathcal{X}}\Lambda\longrightarrow R\Psi^\ad_{\!\mathcal{X}}\Lambda$ by composing
 \begin{align*}
 R\Psi_{\!\mathcal{X}}\Lambda&\yleftarrow[\cong]{\lambda_U}R\Psi_{\!\mathcal{X},U}\Lambda
 =Ra_*i^*R\overline{\jmath}_*\Lambda
 \longrightarrow Ra_*R(\Psi^\ad_{\!\mathcal{X}_{/U_\bullet}})_*\varepsilon^*\Lambda\\
  &=R\Psi^\ad_{\!\mathcal{X}} Ra^\ad_*\Lambda\yleftarrow{\cong}
  R\Psi^\ad_{\!\mathcal{X}}\Lambda.
 \end{align*}
 Here the isomorphy of the last arrow follows from the fact that $a^\ad$ is a morphism of cohomological
 descent.  By the same method as in the proof of Proposition \ref{prop:functoriality}, 
 we can prove that the morphism above is independent of the choice of $U$.

 Since the construction of $\varepsilon^*$ is functorial, in order to
 prove that $\varepsilon^*$ is an isomorphism,
 we may work locally (\cf Proposition \ref{prop:open-base-change}).
 Thus, by Proposition \ref{prop:affine-alg} and Proposition \ref{prop:affine-pi-adic},
 $\varepsilon^*$ is an isomorphism if $\mathcal{X}$ is locally algebraizable or $\varpi$-adic.
\end{prf}

\begin{rem}
 It is plausible that the morphism $\varepsilon^*$ is an isomorphism for a general special formal scheme 
 $\mathcal{X}$ over $\mathcal{S}$ such that $\mathcal{X}_\red$ is separated. 
\end{rem}

\begin{rem}
 One of the most important properties of Berkovich's formal nearby cycle functor is
 the continuity theorem \cite[Theorem 4.1]{MR1395723}. By Theorem \ref{thm:comp-Berkovich},
 our $R\Psi_{\!\mathcal{X}}\Lambda$ also has the property if $\mathcal{X}$ is locally algebraizable or
 $\varpi$-adic. It is an interesting problem to prove the continuity theorem on $R\Psi_{\!\mathcal{X},\mathcal{Z}}\Lambda$ for general $\mathcal{Z}$. It seems difficult, since our functor has no apparent relation to
 rigid geometry.

 In \cite{RZ-GSp4}, we need a continuity of the following type.
 \begin{quote}
  Let $G$ be a locally profinite group which acts continuously on a quasi-compact special formal scheme
  $\mathcal{X}$ over $\mathcal{S}$ such that $\mathcal{X}_\red$ is separated,
  and $\mathcal{Z}=(\mathcal{Z}_1,\mathcal{Z}_2)$ a pair of closed formal subschemes of $\mathcal{X}_s$
  with $\mathcal{Z}_1\subset \mathcal{Z}_2$ which is preserved by $G$.
  Then the action of $G$ on $H^q_c(\mathcal{X}_\red,R\Psi_{\!\mathcal{X},\mathcal{Z}}\Q_\ell)$ is
  smooth (namely, the stabilizer of each element is an open subgroup of $G$).
 \end{quote}
 We do not have a proof of it, but we can prove it when $G$ is obtained as the set of $K$-valued points
 of a linear algebraic group over a $p$-adic field $K$. Our proof is not geometric but purely algebraic;
 we use some properties of pro-$p$ groups.
 See \cite[Section 2]{RZ-GSp4}.
\end{rem}

\subsection{Compactly supported cohomology and $R\Psi_{\!\mathcal{X},c}\Lambda$}
In this subsection, we will relate $H^q_c(\mathcal{X}_\red,R\Psi_{\!\mathcal{X},c}\Lambda)$ 
to the compactly supported cohomology $H^q_c(t(\mathcal{X})_{\overline{\eta}},\Lambda)$ under some condition.
Let $\mathcal{X}$ be a quasi-compact special formal scheme which is separated over $\mathcal{S}$,
namely, the diagonal $\mathcal{X}\longrightarrow \mathcal{X}\times_{\mathcal{S}}\mathcal{X}$ is
a closed immersion.
Then $t(\mathcal{X})_{\eta}$ is separated over $\Spa (F,R)$ (\cf Lemma \ref{lem:generic-fiber} v)).
Moreover, by the following lemma, $t(\mathcal{X})_{\eta}$ is taut:

\begin{lem}\label{lem:special-taut}
 For a quasi-compact special formal scheme $\mathcal{X}$ which is separated over $\mathcal{S}$,
 $t(\mathcal{X})_{\eta}$ and $t(\mathcal{X})_{\overline{\eta}}$ are taut.
\end{lem}

\begin{prf}
 First we will consider the case where $\mathcal{X}=\Spf R[[T]]$. For an integer $m\ge 1$, we put
 \begin{align*}
  \mathbb{U}_m&=\bigl\{v\in t(\mathcal{X})\bigm| v(T^m)\le v(\varpi)\le v(T^{m-1}),\ v(\varpi)\neq 0\bigr\}\\
  &=\{v\in t(\mathcal{X})\bigm| v(T^m)\le v(\varpi)\neq 0\}\cap 
  \{v\in t(\mathcal{X})\bigm| v(\varpi)\le v(T^{m-1})\neq 0\},
 \end{align*}
 which is a quasi-compact open subset of $t(\mathcal{X})_{\eta}$.
 Then it is easy to see that $t(\mathcal{X})_\eta=\bigcup_{m=1}^\infty \mathbb{U}_m$ and 
 $\mathbb{U}_m\cap \mathbb{U}_n=\varnothing$ unless $m-n=\pm 1$ (note that $v(T)<1$ since $v$ is a continuous
 valuation).
 Therefore $t(\mathcal{X})_\eta$ is taut by \cite[Lemma 5.1.3 ii)]{MR1734903}.

 Now we consider the general case. By \cite[Lemma 5.1.3 ii)]{MR1734903},
 we may assume that $\mathcal{X}$ is affine.
 Then $\mathcal{X}$ is a closed formal subscheme of 
 $\Spf R\langle T_1,\ldots,T_m\rangle[[S_1,\ldots,S_n]]$ for some $m$, $n$.
 Since a closed adic subspace of a taut adic space is taut,
 we may assume that $\mathcal{X}=\Spf R\langle T_1,\ldots,T_m\rangle[[S_1,\ldots,S_n]]$.
 Put $\mathcal{Y}=\Spf R[[S_1,\ldots,S_n]]$. Then 
 \[
  t(\mathcal{Y})_\eta=t(\Spf R[[S_1]])_\eta\times_{\Spa(F,R)}\cdots \times_{\Spa(F,R)}t(\Spf R[[S_n]])_\eta
 \]
 is taut by \cite[Lemma 5.1.3 v), Lemma 5.1.4 iii)]{MR1734903}.
 Since $t(\mathcal{X})_\eta\longrightarrow t(\mathcal{Y})_\eta$ is quasi-compact and quasi-separated, 
 $t(\mathcal{X})_\eta$ is also taut \cite[Lemma 5.1.3 iv)]{MR1734903}.
 Thus by \cite[Lemma 5.1.4 iii)]{MR1734903}, $t(\mathcal{X})_{\overline{\eta}}$ is also taut.
\end{prf}

Therefore, we may define the compactly supported cohomology $H^q_c(t(\mathcal{X})_{\overline{\eta}},\Lambda)$
of $t(\mathcal{X})_{\overline{\eta}}$ (\cf \cite[Chapter 5]{MR1734903}).

\begin{cor}\label{cor:H_c-limit}
 We have a natural isomorphism 
 \[
  H^q_c\bigl(t(\mathcal{X})_{\overline{\eta}},\Lambda\bigr)
 =\varinjlim_{(\mathcal{X}'\to \mathcal{X})\in \Phi^{\mathrm{op}}_{\mathcal{X}}}H^q_c\bigl(t(\mathcal{X}'_\varpi)_{\overline{\eta}},\Lambda\bigr).
 \]
\end{cor}

\begin{prf}
 By Lemma \ref{lem:special-taut} and \cite[Lemma 5.1.3]{MR1734903},
 $t(\mathcal{X}'_{\varpi})_{\overline{\eta}}\longrightarrow t(\mathcal{X})_{\overline{\eta}}$ is 
 a taut morphism between taut adic spaces.
 Therefore the corollary follows from \cite[Proposition 5.4.5]{MR1734903}.
\end{prf}

\subsubsection{Complements on the functor $Rf_!$}
For a separated morphism $f\colon X\longrightarrow Y$ locally of $^+$weakly finite type between 
analytic adic spaces, Huber constructed the functor $Rf_!$ \cite[Section 5]{MR1734903} (in \cite{MR1734903},
it is denoted by $R^+f_!$). If $f$ is partially proper or an immersion,
then $Rf_!$ is the derived functor of $f_!$ (in the latter case we moreover have $Rf_!=f_!$).
In these cases, we have a natural morphism of functors
$Rf_!\longrightarrow Rf_*$, for $f_!$ is a subfunctor of $f_*$.

We need the following proposition:

\begin{prop}\label{prop:cpt-to-ord}
 For a separated morphism $f\colon X\longrightarrow Y$ locally of $^+$weakly finite type between 
 analytic adic spaces, we have a morphism of functors $Rf_!\longrightarrow Rf_*$ which satisfies
 the following conditions:
 \begin{itemize}
  \item If $f$ is partially proper or an immersion, then $Rf_!\longrightarrow Rf_*$ coincides with
	the morphism above.
  \item The morphism is compatible with composition. Namely, we have the following commutative diagram:
	\[
	 \xymatrix{%
	R(g\circ f)_!\ar@{<->}[r]^{\cong}\ar[d]& Rg_!\circ Rf_!\ar[d]\\
	R(g\circ f)_*\ar@{=}[r]& Rg_*\circ Rf_*\lefteqn{.}
	}
	\]
 \end{itemize}
\end{prop}

Note that the natural morphism of functors $f_!\longrightarrow f_*$ is compatible with composition
(\cf \cite[Proposition 5.2.2 iii)]{MR1734903}). Therefore, if both $f$ and $g$ are partially proper
or are immersions, then the second condition in the proposition above is satisfied.
Thus, as in the proof of \cite[Theorem 5.4.3]{MR1734903}, it suffices to show the following lemma:

\begin{lem}
 Consider the following commutative diagram of analytic adic spaces where $f$, $f'$ are
 partially proper and $j$, $j'$ are immersions:
 \[
  \xymatrix{%
 Y'\ar[r]^-{j'}\ar[d]^-{f'}& Y\ar[d]^-{f}\\
 X'\ar[r]^-{j}& X\lefteqn{.}
 }
 \]
 Then two morphisms of functors $j_!\circ Rf'_!\longrightarrow Rj_*\circ Rf'_*$ and 
 $Rf_!\circ j'_!\longrightarrow Rf_*\circ Rj'_*$ are equal under the identifications
 $j_!\circ Rf'_!=Rf_!\circ j'_!$ in \cite[Lemma 5.4.2]{MR1734903} and
 $Rj_*\circ Rf'_*=Rf_*\circ Rj'_*$.
\end{lem}

\begin{prf}
 Obviously $j_!\circ Rf'_!$ is the derived functor of $(j\circ f')_!$.
 By \cite[Lemma 5.4.2]{MR1734903}, $Rf_!\circ j'$ is the derived functor of $(f\circ j')_!$.
 Now it is easy see that the two morphisms in the lemma
 can be identified with the morphism induced by $g_!\longrightarrow g_*$, where
 $g=j\circ f'=f\circ j'$.
\end{prf}

\subsubsection{Pseudo-compactifications}
In order to construct the comparison homomorphism
\[
 \varepsilon_*\colon H^q_c\bigl(t(\mathcal{X})_{\overline{\eta}},\Lambda\bigr)\longrightarrow 
 H^q_c(\mathcal{X}_\red,R\Psi_{\!\mathcal{X},c}\Lambda),
\]
we introduce the notion of pseudo-compactifications.
We will assume that every formal scheme appearing here is a quasi-compact special formal scheme
which is separated over $\mathcal{S}$.

\begin{defn}\label{defn:pseudo-proper}
 We say that a formal scheme $\mathcal{X}$ is pseudo-proper over $\mathcal{S}$ 
 if $\mathcal{X}_\red$ is proper over $\Spec k$ and $t(\mathcal{X})_\eta$ is partially proper over $\Spa (F,R)$. 
\end{defn}

\begin{lem}\label{lem:pseudo-proper}
 \begin{enumerate}
  \item Let $X$ be a proper $S$-scheme and $Y$ a closed subscheme of $X_s$. Then the completion $\mathcal{X}$
	of $X$ along $Y$ is pseudo-proper over $\mathcal{S}$. 
  \item Let $\mathcal{X}$ be a formal scheme which is pseudo-proper over $\mathcal{S}$.
	For an object $\mathcal{X}'\longrightarrow \mathcal{X}$ of $\Phi_\mathcal{X}$,
	$\mathcal{X}'$ is pseudo-proper over $\mathcal{S}$.
  \item For formal schemes $\mathcal{X}$ and $\mathcal{Y}$ which are pseudo-proper over $\mathcal{S}$,
	$\mathcal{X}\times_{\mathcal{S}}\mathcal{Y}$ is also pseudo-proper over $\mathcal{S}$.
 \end{enumerate}
\end{lem}

\begin{prf}
 For i), we have only to verify that $t(\mathcal{X})_\eta$ is partially proper over $\Spa (F,R)$.
 It is equivalent to that for every quasi-compact open subset $\mathbb{U}$ of $t(\mathcal{X})_\eta$,
 the pseudo-adic space $(t(\mathcal{X})_\eta,\overline{\mathbb{U}})$ is proper over $\Spa (F,R)$,
 where $\overline{\mathbb{U}}$ denotes the closure of $\mathbb{U}$ in $t(\mathcal{X})_\eta$.
 Let us denote the $\varpi$-adic completion of $X$ by $\mathcal{X}_1$. Then by \cite[Lemma 3.13]{MR1620118},
 $t(\mathcal{X})_\eta$ is the interior of a closed constructible subset of $t(\mathcal{X}_1)_\eta$.
 Therefore $\overline{\mathbb{U}}$ coincides with the closure of $\mathbb{U}$ in $t(\mathcal{X}_1)_\eta$
 (\cf \cite[proof of Lemma 1.3 iii)]{MR1620114}). Since $t(\mathcal{X}_1)_\eta$ is proper over $\Spa (F,R)$, 
 $(t(\mathcal{X}_1)_\eta,\overline{\mathbb{U}})$ is also proper over $\Spa(F,R)$, and thus 
 $(t(\mathcal{X})_\eta,\overline{\mathbb{U}})$ is also proper over $\Spa(F,R)$.

 ii) is clear from $t(\mathcal{X}')_\eta=t(\mathcal{X})_\eta$.
 iii) is a consequence of Lemma \ref{lem:generic-fiber} v).
\end{prf}

\begin{rem}
 It is plausible that a formal scheme $\mathcal{X}$ over $\mathcal{S}$ is pseudo-proper
 if $\mathcal{X}_\red$ is proper over $\Spec k$.
\end{rem}

\begin{defn}\label{defn:pseudo-comp}
 \begin{enumerate}
  \item A pseudo-compactification of $\mathcal{X}$ is an immersion 
	$\mathcal{X}\hooklongrightarrow \overline{\mathcal{X}}$ into a formal scheme $\overline{\mathcal{X}}$
	which is pseudo-proper over $\mathcal{S}$. 
	A formal scheme which has a pseudo-compactification is said to be pseudo-compactifiable.
  \item An adic $\mathcal{S}$-morphism of formal schemes $f\colon \mathcal{Y}\longrightarrow \mathcal{X}$ is
	said to be pseudo-compactifiable if there exists a diagram of formal schemes over $\mathcal{S}$
	\[
	\xymatrix{%
	\mathcal{Y}\ar@{^(->}[r]\ar[d]^-{f}& \overline{\mathcal{Y}}\ar[d]^-{f'}\\
	\mathcal{X}\ar@{^(->}[r]& \overline{\mathcal{X}}\lefteqn{,}
	}
	\]
	where the horizontal arrows are pseudo-compactifications and $f'$ is adic (hence proper,
	since $f'_\red$ is proper).
 \end{enumerate}
\end{defn}

\begin{exa}\label{exa:pseudo-comp}
 \begin{enumerate}
  \item If $\mathcal{X}$ is affine, then it is pseudo-compactifiable.
	Indeed, $\mathcal{X}$ is a closed formal subscheme of 
	$\Spf R\langle T_1,\ldots,T_m\rangle[[S_1,\ldots,S_n]]$ for some $m$, $n$,
	which is an open formal subscheme of 
	$\widehat{\P}^m_{\mathcal{S}}\times_{\mathcal{S}}\Spf R[[S_1,\ldots,S_n]]$. The last formal scheme
	is pseudo-proper over $\mathcal{S}$ by Lemma \ref{lem:pseudo-proper} i), since it
	is the completion of $\P^m_S\times_S \P^n_S$ along $\P^m_k\times \{0\}$.
  \item If $\mathcal{X}$ is algebraizable, then it is pseudo-compactifiable.
	Indeed, for an $S$-scheme $X$ which is separated of finite type and
	a closed subscheme $Y$ of $X_s$ such that the completion of $X$ along $Y$ is isomorphic to $\mathcal{X}$,
	take a compactification $\overline{X}$ of $X$ and the closure $\overline{Y}$ of $Y$ in $\overline{X}$.
	Then the completion $\overline{\mathcal{X}}$ of $\overline{X}$ along $\overline{Y}$ is pseudo-proper
	over $\mathcal{S}$ (Lemma \ref{lem:pseudo-proper} i))
	and contains $\mathcal{X}$ as an open formal subscheme.
  \item If $\mathcal{X}$ is pseudo-compactifiable, so are its open formal subschemes.
  \item If $\mathcal{X}$ is pseudo-compactifiable, then an isomorphism
	of formal schemes $\mathcal{Y}\yrightarrow{\cong}\mathcal{X}$ is
	pseudo-compactifiable. Indeed, for a pseudo-compactification
	$\mathcal{X}\hooklongrightarrow\overline{\mathcal{X}}$ of $\mathcal{X}$, 
	$\mathcal{Y}\yrightarrow{\cong}\mathcal{X}\hooklongrightarrow\overline{\mathcal{X}}$ gives
	a pseudo-compactification of $\mathcal{Y}$.
	More generally, an immersion $\mathcal{Y}\longrightarrow\mathcal{X}$ is
	pseudo-compactifiable.
 \end{enumerate}
\end{exa}

\begin{lem}\label{lem:pseudo-comp-adm-bu}
 \begin{enumerate}
  \item Let $\mathcal{X}$ be a pseudo-compactifiable formal scheme over $\mathcal{S}$.
	Then, for every object $\mathcal{X}'\longrightarrow \mathcal{X}$ of $\Phi_\mathcal{X}$,
	$\mathcal{X}'$ is also pseudo-compactifiable. Moreover, 
	for every morphism $(\mathcal{X}''\longrightarrow \mathcal{X})\longrightarrow (\mathcal{X}'\longrightarrow \mathcal{X})$ of $\Phi_\mathcal{X}$, the morphism $\mathcal{X}''\longrightarrow \mathcal{X}'$ is
	pseudo-compactifiable.
  \item Let $f\colon \mathcal{Y}\longrightarrow \mathcal{X}$ be a pseudo-compactifiable morphism of
	pseudo-compactifiable formal schemes over $\mathcal{S}$. Then, 
	for every object $\mathcal{X}'\longrightarrow \mathcal{X}$ of $\Phi_\mathcal{X}$,
	there exist an object $\mathcal{Y}'\longrightarrow \mathcal{Y}$ of $\Phi_\mathcal{Y}$ and
	a pseudo-compactifiable morphism $f'\colon \mathcal{Y}'\longrightarrow \mathcal{X}'$
	such that the following diagram is commutative:
	\[
	\xymatrix{%
	\mathcal{Y}'\ar[r]^-{f'}\ar[d]& \mathcal{X'}\ar[d]\\
	\mathcal{Y}\ar[r]^-{f}& \mathcal{X}\lefteqn{.}
	}
	\]
 \end{enumerate}
\end{lem}

\begin{prf}
 We prove i). By definition, there exist a closed immersion $\mathcal{X}\hooklongrightarrow \mathcal{X}_1$
 and an open immersion $\mathcal{X}_1\hooklongrightarrow \overline{\mathcal{X}}$
 such that $\overline{\mathcal{X}}$ is pseudo-proper over $\mathcal{S}$.
 We denote the ideal of the center of the blow-up $\mathcal{X}'\longrightarrow \mathcal{X}$ by $\mathcal{J}$.
 Let $\mathcal{J}_1$ be the ideal of $\mathcal{O}_{\mathcal{X}_1}$ that is naturally induced from $\mathcal{J}$
 and $\mathcal{X}'_1\longrightarrow \mathcal{X}_1$ the admissible blow-up along $\mathcal{J}_1$.
 Then, we have a natural closed immersion $\mathcal{X}'\hooklongrightarrow \mathcal{X}_1'$.
 Moreover, by \cite[I, (9.4.7)]{EGA}, there exists an extension of $\mathcal{J}_1$ to an admissible
 ideal of $\overline{\mathcal{X}}$, which we denote by $\overline{\mathcal{J}}$.
 If we denote the admissible blow-up along $\overline{\mathcal{J}}$ by $\overline{\mathcal{X}}'\longrightarrow \overline{\mathcal{X}}$, we have the following commutative diagram:
 \[
  \xymatrix{%
 \mathcal{X}'\ar[r]\ar[d]& \mathcal{X}'_1\ar[r]\ar[d]& \overline{\mathcal{X}}'\ar[d]\\
 \mathcal{X}\ar[r]& \mathcal{X}_1\ar[r]& \overline{\mathcal{X}}\lefteqn{.}
 }
 \]
 The immersion $\mathcal{X}'\hooklongrightarrow \overline{\mathcal{X}}'$ gives a pseudo-compactification of
 $\mathcal{X}'$ by Lemma \ref{lem:pseudo-proper} ii).
 Moreover, the diagram above shows that the admissible blow-up
 $\mathcal{X}'\longrightarrow \mathcal{X}$ is pseudo-compactifiable.
 Denote the ideal of the center of the blow-up $\mathcal{X}''\longrightarrow \mathcal{X}$ by $\mathcal{K}$.
 Then $\mathcal{X}''\longrightarrow \mathcal{X}'$ is the admissible blow-up of $\mathcal{X}'$
 along $\mathcal{K}\mathcal{O}_{\mathcal{X}'}$. Hence it is compactifiable by the argument above.

 For ii), we simply take $\mathcal{Y}'\longrightarrow \mathcal{Y}$ as the admissible blow-up 
 along $(f^{-1}\mathcal{J})\mathcal{O}_{\mathcal{Y}}$.
\end{prf}

\begin{defn}\label{def:cosp}
 Let $\mathcal{X}$ be a pseudo-compactifiable $\varpi$-adic formal scheme over $\mathcal{S}$ and
 $L$ an object of $D^+(t(\mathcal{X})_{\overline{\eta}},\Lambda)$.
 We take a pseudo-compactification $\mathcal{X}\stackrel{j}{\hooklongrightarrow}\overline{\mathcal{X}}$ and
 define the homomorphism 
 \[
 \xi_{\mathcal{X}}\colon H^q_c(\mathcal{X}_\red,R\Psi^\ad_{\!\mathcal{X}}L)\longrightarrow H^q_c\bigl(t(\mathcal{X})_{\overline{\eta}},L\bigr)
 \]
as the composite of
 \begin{align*}
  H^q_c(\mathcal{X}_\red,R\Psi^\ad_{\!\mathcal{X}} L)
  &=H^q(\overline{\mathcal{X}}_\red,j_{\red!}R\Psi^\ad_{\!\mathcal{X}}L)
  \yrightarrow{(1)} H^q\bigl(\overline{\mathcal{X}}_\red,R\Psi^\ad_{\!\overline{\mathcal{X}}} t(j)_{\overline{\eta}!}L\bigr)\\
  &=H^q\bigl(t(\overline{\mathcal{X}})_{\overline{\eta}},t(j)_{\overline{\eta}!}L\bigr)
  \yleftarrow[\cong]{(2)}H^q_c\bigl(t(\mathcal{\overline{X}})_{\overline{\eta}},t(j)_{\overline{\eta}!}L\bigr)
  =H^q_c\bigl(t(\mathcal{X})_{\overline{\eta}},L\bigr).
 \end{align*}
 Here (1) is induced from the natural isomorphism
 $R\Psi^\ad_{\!\mathcal{X}}L\cong j_\red^*R\Psi^\ad_{\!\overline{\mathcal{X}}} t(j)_{\overline{\eta}!}L$.
 The homomorphism $(2)$ is the one constructed in Proposition \ref{prop:cpt-to-ord}.
 It is an isomorphism since $t(\mathcal{X})_\eta$ is a quasi-compact open subset of
 $t(\overline{\mathcal{X}})_\eta$ (recall that $\mathcal{X}$ is $\varpi$-adic)
 and thus the closure of $t(\mathcal{X})_\eta$ in $t(\overline{\mathcal{X}})_\eta$
 is proper over $\Spa(F,R)$ as a pseudo-adic space
 (\cf proof of Lemma \ref{lem:pseudo-proper} i)).
\end{defn}

\begin{lem}\label{lem:xi-property}
 Let $\mathcal{X}$ and $L$ be as in Definition \ref{def:cosp}.
 \begin{enumerate}
 \item The homomorphism $\xi_{\mathcal{X}}$ is independent of the choice of a pseudo-compactification.
  \item For a pseudo-compactifiable proper morphism $f\colon \mathcal{Y}\longrightarrow \mathcal{X}$ between
 pseudo-compactifiable $\varpi$-adic formal schemes over $\mathcal{S}$, the following diagram is commutative:
	\[
 \xymatrix{%
 H^q_c(\mathcal{X}_\red,R\Psi^\ad_{\!\mathcal{X}} L)\ar[d]^-{f^*}\ar[r]^-{\xi_{\mathcal{X}}}& H^q_c(t(\mathcal{X})_{\overline{\eta}},L)\ar[d]^-{t(f)_{\overline{\eta}}^*}\\
 H^q_c\bigl(\mathcal{Y}_\red,R\Psi^\ad_{\!\mathcal{Y}} t(f)_{\overline{\eta}}^*L\bigr)\ar[r]^-{\xi_{\mathcal{Y}}}& H^q_c\bigl(t(\mathcal{Y})_{\overline{\eta}},t(f)_{\overline{\eta}}^*L\bigr)\lefteqn{.}
	}
	\]
  \item For an open immersion $\mathcal{U}\hooklongrightarrow \mathcal{X}$, the following diagram is commutative:
	\[
	\xymatrix{%
	H^q_c\bigl(\mathcal{U}_\red,R\Psi^\ad_{\mathcal{U}} (L\vert_{t(\mathcal{U})_{\overline{\eta}}})\bigr)\ar[d]\ar[r]^-{\xi_{\mathcal{U}}}& H^q_c\bigl(t(\mathcal{U})_{\overline{\eta}},L\vert_{t(\mathcal{U})_{\overline{\eta}}}\bigr)\ar[d]\\
	H^q_c(\mathcal{X}_\red,R\Psi^\ad_{\!\mathcal{X}} L)\ar[r]^-{\xi_{\mathcal{X}}}& H^q_c\bigl(t(\mathcal{X})_{\overline{\eta}},L\bigr)\lefteqn{.}
	}
	\]
  \item For an open covering $\mathcal{X}=\bigcup_{i\in I}\mathcal{U}_i$ of $\mathcal{X}$,
	we have the following morphism of spectral sequences
	(\cf \cite[Remark 5.5.12 iii)]{MR1734903}), where we put $\mathcal{U}_{i_0,\ldots,i_p}=\mathcal{U}_{i_0}\cap \cdots \cap \mathcal{U}_{i_p}$:
	\[
  \xymatrix{
 E_1^{-p,q}=\bigoplus_{(i_0,\ldots,i_p)\in I^{p+1}}H^q_c\bigl((\mathcal{U}_{i_0,\ldots,i_p})_\red,R\Psi^\ad_{\mathcal{U}_{i_0,\ldots,i_p}} L\bigr)\ar@{=>}[r]\ar[d]^-{\oplus\xi_{\mathcal{U}_{i_0,\ldots,i_p}}}& H^{-p+q}_c(\mathcal{X}_\red,R\Psi^\ad_{\!\mathcal{X}} L)\ar[d]^-{\xi_{\mathcal{X}}}\\
 E_1^{-p,q}=\bigoplus_{(i_0,\ldots,i_p)\in I^{p+1}}H^q_c\bigl(t(\mathcal{U}_{i_0,\ldots,i_p})_{\overline{\eta}},L\bigr)\ar@{=>}[r]& 
 H^{-p+q}_c\bigl(t(\mathcal{X})_{\overline{\eta}},L\bigr)\lefteqn{.}
 }
 \]
 \end{enumerate}
\end{lem}

\begin{prf}
 We will prove i). Let $j'\colon \mathcal{X}\hooklongrightarrow \overline{\mathcal{X}}'$ be another pseudo-compactification.
 Then the composite
 $\mathcal{X}\hooklongrightarrow \mathcal{X}\times_{\mathcal{S}}\mathcal{X}\yrightarrow{j\times j'}\overline{\mathcal{X}}\times_{\mathcal{S}}\overline{\mathcal{X}}'$ is also a pseudo-compactification by 
 Lemma \ref{lem:pseudo-proper} iii).
 Thus we may assume that there exists a morphism $\pi\colon \overline{\mathcal{X}}'\longrightarrow \overline{\mathcal{X}}$ over $\mathcal{S}$ such that $\pi\circ j'=j$.
 Obviously $\pi_\red$ is proper.
 The claim is an immediate consequence of the commutativity of
 the following diagrams (we also write $j$ and $j'$ for $t(j)_{\overline{\eta}}$ and
 $t(j')_{\overline{\eta}}$):
 \[
 \xymatrix{%
 H^q_c\bigl(t(\mathcal{X})_{\overline{\eta}},L\bigr)\ar@{=}[r]\ar@{=}[ddd]
 &H^q_c\bigl(t(\overline{\mathcal{X}})_{\overline{\eta}},j_!L\bigr)\ar[r]^-{\cong}\ar@{=}[d]
 &H^q\bigl(t(\overline{\mathcal{X}})_{\overline{\eta}},j_!L\bigr)\ar@{=}[d]\\
 &H^q_c\bigl(t(\overline{\mathcal{X}})_{\overline{\eta}},R\pi_!j'_!L\bigr)\ar[r]\ar@{=}[dd]
 &H^q\bigl(t(\overline{\mathcal{X}})_{\overline{\eta}},R\pi_!j'_!L\bigr)\ar[d]\\
 &&H^q\bigl(t(\overline{\mathcal{X}})_{\overline{\eta}},R\pi_*j'_!L\bigr)\ar@{=}[d]\\
 H^q_c\bigl(t(\mathcal{X})_{\overline{\eta}}, L\bigr)\ar@{=}[r]
 &H^q_c\bigl(t(\overline{\mathcal{X}}')_{\overline{\eta}},j'_!L\bigr)\ar[r]^-{\cong}
 &H^q\bigl(t(\overline{\mathcal{X}}')_{\overline{\eta}},j'_!L\bigr)\lefteqn{,}
 }
 \]
 \[
 \xymatrix{%
 H^q(\overline{\mathcal{X}}_\red,j_{\red!}R\Psi^\ad_{\!\mathcal{X}}L)\ar[r]\ar@{=}[dd]
 &H^q(\overline{\mathcal{X}}_\red,R\Psi^\ad_{\!\overline{\mathcal{X}}}j_!L)\ar[d]\ar@{=}[r]
 &H^q\bigl(t(\overline{\mathcal{X}})_{\overline{\eta}},j_!L\bigr)\ar[d]\\
 &H^q(\overline{\mathcal{X}}_\red,R\Psi^\ad_{\!\overline{\mathcal{X}}} R\pi_* j'_!L)\ar@{=}[r]\ar@{=}[d]
 &H^q\bigl(t(\overline{\mathcal{X}})_{\overline{\eta}},R\pi_*j'_!L\bigr)\ar@{=}[dd]\\
 H^q(\overline{\mathcal{X}}_\red,R\pi_{\red*}j'_{\red!}R\Psi^\ad_{\!\mathcal{X}}L)\ar[r]\ar@{=}[d]
 &H^q(\overline{\mathcal{X}}_\red,R\pi_{\red*}R\Psi^\ad_{\!\overline{\mathcal{X}}'}j'_!L)\ar@{=}[d]\\
 H^q(\overline{\mathcal{X}}'_\red,j'_{\red!}R\Psi^\ad_{\!\mathcal{X}}L)\ar[r]
 &H^q(\overline{\mathcal{X}}'_\red,R\Psi^\ad_{\!\overline{\mathcal{X}}'} j'_!L) \ar@{=}[r]
 &H^q\bigl(t(\overline{\mathcal{X}}')_{\overline{\eta}},j'_!L\bigr)\lefteqn{.}
 }
 \]
 The first diagram is clearly commutative.
 The commutativity of the second diagram is also trivial, except for the upper left rectangle.
 Let us observe the commutativity of the remained rectangle.
 Take a factorization $\mathcal{X}\stackrel{i}{\hooklongrightarrow}\mathcal{X}_1\stackrel{j_1}{\hooklongrightarrow} \overline{\mathcal{X}}$ of $j$, where $i$ is a closed immersion and $j_1$ is an open immersion.
 Consider the following diagram, whose rectangles are cartesian:
 \[
  \xymatrix{%
 \mathcal{X}\ar[r]^-{\sigma}\ar@{=}[rd]& \mathcal{X}'\ar[r]^-{i'}\ar[d]^-{\pi'}& \mathcal{X}'_1\ar[r]^-{j_1'}\ar[d]^-{\pi_1}&
 \overline{\mathcal{X}}'\ar[d]^-{\pi}\\
 & \mathcal{X}\ar[r]^-{i}& \mathcal{X}_1\ar[r]^-{j_1}&\overline{\mathcal{X}}\lefteqn{.}
 }
 \]
 Note that $\sigma$ is a closed immersion, since $\pi'$ is separated.
 By the adjointness of $(j_1)_{\red!}$ and $(j_1)_\red^*$,
 it suffices to show that the following diagram is commutative
 (note that $i_!=i_*$, $i'_!=i'_*$ and $\sigma_!=\sigma_*$, since $i$, $i'$, $\sigma$ are closed immersions):
 \[
 \xymatrix{%
 i_{\red!}R\Psi^\ad_{\!\mathcal{X}} L\ar@{=}[r]\ar@{=}[ddd]& R\Psi^\ad_{\!\mathcal{X}_1}i_!L\ar@{=}[d]^-{(1)} \\
 &R\Psi^\ad_{\!\mathcal{X}_1} R\pi_{1!}i'_!\sigma_!L\ar[d]^-{(2)}\\
 &R\Psi^\ad_{\!\mathcal{X}_1} R\pi_{1*}i'_*\sigma_*L\ar@{=}[d]\\
 R(\pi_1)_{\red*}i'_{\red*}\sigma_{\red*}R\Psi^\ad_{\!\mathcal{X}} L\ar[r]
 &R(\pi_1)_{\red*} R\Psi^\ad_{\!\mathcal{X}'_1} i'_*\sigma_*L \lefteqn{.}
 }
 \]
 By Proposition \ref{prop:cpt-to-ord}, the composite of (1) and (2) is the evident map
 $R\Psi^\ad_{\!\mathcal{X}_1} i_*L\yrightarrow{=} R\Psi^\ad_{\!\mathcal{X}_1} R\pi_{1*}i'_*\sigma_*L$
 induced from $i=\pi_1\circ i'\circ \sigma$.
 Now the commutativity is clear, since two morphisms we should compare are both obtained as the composites
 of evident identifications such as $g_*\circ f_*=(g\circ f)_*$, where $f$ and $g$ are morphisms of sites.

 The proofs of ii) and iii) are straightforward, and iv) is an easy consequence of
 \cite[Corollary 3.5.11 ii)]{MR1734903}. We will omit them.
\end{prf}

\begin{cor}\label{cor:xi-isom}
 Let $\mathcal{X}$ and $L$ be as in Definition \ref{def:cosp}. Then $\xi_\mathcal{X}$ is an isomorphism.
\end{cor}

\begin{prf}
 By Lemma \ref{lem:xi-property} iv), we may assume that $\mathcal{X}$ is affine.
 Then the isomorphy of $\xi_\mathcal{X}$ follows from \cite[Lemma 2.13]{MR1620114}.
\end{prf}

\begin{rem}\label{rem:xi-isom}
 If $\mathcal{X}$ is not $\varpi$-adic, 
 $H^q_c(t(\mathcal{X})_{\overline{\eta}},L)$ and $H^q_c(\mathcal{X}_\red,R\Psi^\ad_{\!\mathcal{X}}L)$ are
 not isomorphic in general.
 For example, if $\mathcal{X}=\Spf R[[T]]$, then $H^2_c(t(\mathcal{X})_{\overline{\eta}},\Lambda)\cong\Lambda(-1)\neq 0$, while $H^2_c(\mathcal{X}_\red,R\Psi^\ad_{\!\mathcal{X}}\Lambda)=H^2_c(\mathcal{X}_\red,\Lambda)=0$.
\end{rem}

\subsubsection{Construction of the comparison map}
Now we can give the definition of
$\varepsilon_*\colon H^q_c(t(\mathcal{X})_{\overline{\eta}},\Lambda)\longrightarrow H^q_c(\mathcal{X}_\red,R\Psi_{\!\mathcal{X},c}\Lambda)$.
In the sequel, let $\mathcal{X}$ be a quasi-compact special formal scheme
which is separated over $\mathcal{S}$.

\begin{defn}\label{defn:e-X}
 Assume that $\mathcal{X}$ is pseudo-compactifiable. 
 For an object $\mathcal{X}'\longrightarrow \mathcal{X}$ of $\Phi_{\mathcal{X}}$, we define 
 the homomorphism 
 \[
  \varepsilon_{\mathcal{X}'}\colon H^q_c\bigl(t(\mathcal{X}'_\varpi)_{\overline{\eta}},\Lambda\bigr)\longrightarrow H^q_c(\mathcal{X}_\red,R\Psi_{\!\mathcal{X},c}\Lambda)
 \]
 as the composite
 \begin{align*}
  H^q_c\bigl(t(\mathcal{X}'_\varpi)_{\overline{\eta}},\Lambda\bigr)&\yrightarrow[\cong]{\xi^{-1}_{\mathcal{X}'_\varpi}}
  H^q_c\bigl((\mathcal{X}'_\varpi)_\red,R\Psi^\ad_{\!\mathcal{X}'_\varpi}\Lambda\bigr)
  \stackrel{(1)}{\cong}H^q_c\bigl((\mathcal{X}'_\varpi)_\red,R\Psi_{\!\mathcal{X}'_{\varpi}}\Lambda\bigr)\\
  &\yrightarrow{(2)} H^q_c(\mathcal{X}'_\red,R\Psi_{\!\mathcal{X}',c}\Lambda)
  \yleftarrow[\cong]{(3)}H^q_c(\mathcal{X}_\red,R\Psi_{\!\mathcal{X},c}\Lambda).
 \end{align*}
 Here (1) is the isomorphism in Theorem \ref{thm:comp-Berkovich}.
 The morphism (2) is obtained by Proposition \ref{prop:open-base-change};
 note that $R\Psi_{\!\mathcal{X}'_{\varpi},c}\Lambda=R\Psi_{\!\mathcal{X}'_{\varpi}}\Lambda$ since
 $\mathcal{X}'_{\varpi}$ is $\varpi$-adic. The isomorphism (3) is
 due to Proposition \ref{prop:inv-adm-bu}.
 Finally, notice that $\mathcal{X}'_{\varpi}$ is pseudo-compactifiable by Lemma \ref{lem:pseudo-comp-adm-bu} i) and
 Example \ref{exa:pseudo-comp}, and thus $\xi_{\mathcal{X}'_\varpi}$ is defined.
\end{defn}

\begin{lem}
 Assume that $\mathcal{X}$ is pseudo-compactifiable.
 Let $(\mathcal{X}''\longrightarrow \mathcal{X})\longrightarrow (\mathcal{X}'\longrightarrow \mathcal{X})$
 be a morphism in $\Phi_{\mathcal{X}}$. Then we have the following commutative diagram:
 \[
  \xymatrix{%
 H^q_c\bigl(t(\mathcal{X}'_\varpi)_{\overline{\eta}},\Lambda\bigr)\ar[r]^-{\varepsilon_{\mathcal{X}'}}\ar[d]
 & H^q_c(\mathcal{X}_\red,R\Psi_{\!\mathcal{X},c}\Lambda)\ar@{=}[d]\\
 H^q_c\bigl(t(\mathcal{X}''_\varpi)_{\overline{\eta}},\Lambda\bigr)\ar[r]^-{\varepsilon_{\mathcal{X}''}}
 & H^q_c(\mathcal{X}_\red,R\Psi_{\!\mathcal{X},c}\Lambda)\lefteqn{.}
 }
 \]
\end{lem}

\begin{prf}
 Let $\mathcal{U}''$ be the inverse image of $\mathcal{X}'_\varpi$ by 
 $\mathcal{X}''\longrightarrow \mathcal{X}'$. Consider the following diagrams:
 \[
  \xymatrix{%
 H^q_c\bigl(t(\mathcal{X}'_\varpi)_{\overline{\eta}},\Lambda\bigr)\ar[d]&
 H^q_c\bigl((\mathcal{X}'_\varpi)_\red,R\Psi^\ad_{\!\mathcal{X}'_\varpi}\Lambda\bigr)\ar[l]_-{\xi_{\mathcal{X}'_{\varpi}}}\ar[d]&
 H^q_c\bigl((\mathcal{X}'_\varpi)_\red,R\Psi_{\!\mathcal{X}'_\varpi}\Lambda\bigr)\ar[l]_-{\cong}\ar[d]\\
 H^q_c\bigl(t(\mathcal{U}'')_{\overline{\eta}},\Lambda\bigr)\ar[d]
 & H^q_c\bigl((\mathcal{U}'')_\red,R\Psi^\ad_{\mathcal{U}''}\Lambda\bigr)\ar[l]_-{\xi_{\mathcal{U}''}}\ar[d]
 & H^q_c\bigl((\mathcal{U}'')_\red,R\Psi_{\mathcal{U}''}\Lambda\bigr)\ar[l]_-{\cong}\ar[d]\\
 H^q_c\bigl(t(\mathcal{X}''_\varpi)_{\overline{\eta}},\Lambda\bigr)&
 H^q_c\bigl((\mathcal{X''_\varpi})_\red,R\Psi^\ad_{\!\mathcal{X}''_\varpi}\Lambda\bigr)\ar[l]_-{\xi_{\mathcal{X}''_{\varpi}}}&
 H^q_c\bigl((\mathcal{X}''_\varpi)_\red,R\Psi_{\mathcal{X}''_\varpi}\Lambda\bigr)\ar[l]_-{\cong}\lefteqn{,}
 }
 \]
 \[
  \xymatrix{%
 H^q_c\bigl((\mathcal{X}'_\varpi)_\red,R\Psi_{\!\mathcal{X}'_\varpi}\Lambda\bigr)\ar[r]\ar[d]
 & H^q_c(\mathcal{X}'_\red,R\Psi_{\!\mathcal{X}',c}\Lambda)\ar[dd]
 & H^q_c(\mathcal{X}_\red,R\Psi_{\!\mathcal{X},c}\Lambda)\ar[l]_-{\cong}\ar@{=}[dd]\\
 H^q_c\bigl((\mathcal{U}'')_\red,R\Psi_{\mathcal{U}''}\Lambda\bigr)\ar[d]\\
 H^q_c\bigl((\mathcal{X}''_\varpi)_\red,R\Psi_{\!\mathcal{X}''_\varpi}\Lambda\bigr)\ar[r]&
 H^q_c(\mathcal{X}''_\red,R\Psi_{\!\mathcal{X}'',c}\Lambda)&
 H^q_c(\mathcal{X}_\red,R\Psi_{\!\mathcal{X},c}\Lambda)\ar[l]_-{\cong}\lefteqn{.}
 }
 \]
 The upper diagram is commutative by Lemma \ref{lem:xi-property} ii), iii); note that
 $\mathcal{U}''\longrightarrow \mathcal{X}'_\varpi$ is pseudo-compactifiable 
 by Lemma \ref{lem:pseudo-comp-adm-bu} i).
 The commutativity of the lower diagram follows from easy consideration.
\end{prf}

\begin{defn}
 Assume that $\mathcal{X}$ is pseudo-compactifiable.
 Then we define 
 \[
 \varepsilon_*\colon H^q_c\bigl(t(\mathcal{X})_{\overline{\eta}},\Lambda\bigr)\longrightarrow H^q_c(\mathcal{X}_\red,R\Psi_{\!\mathcal{X},c}\Lambda)
 \]
 as the composite of
 \[
  H^q_c\bigl(t(\mathcal{X})_{\overline{\eta}},\Lambda\bigr)\yleftarrow{\cong}\varinjlim_{(\mathcal{X}'\to \mathcal{X})\in \Phi^{\mathrm{op}}_{\mathcal{X}}}H^q_c\bigl(t(\mathcal{X}'_\varpi)_{\overline{\eta}},\Lambda\bigr)
 \yrightarrow{\varinjlim \varepsilon_{\mathcal{X}'}}H^q_c(\mathcal{X}_\red,R\Psi_{\!\mathcal{X},c}\Lambda).
 \]
\end{defn}

It is easy to see that $\varepsilon_*$ is functorial with respect to the push-forward by an open immersion.
We can also prove the functoriality for the pull-back by a pseudo-compactifiable proper morphism 
of formal schemes $\mathcal{Y}\longrightarrow \mathcal{X}$:

\begin{prop}\label{prop:e-funct-pull}
 Let $f\colon \mathcal{Y}\longrightarrow \mathcal{X}$ be a pseudo-compactifiable proper morphism
 of pseudo-compactifiable formal schemes over $\mathcal{S}$. 
 Then, we have the following commutative diagram:
 \[
  \xymatrix{%
 H^q_c\bigl(t(\mathcal{X})_{\overline{\eta}},\Lambda\bigr)\ar[r]^-{\varepsilon_*}\ar[d]^-{t(f)_{\overline{\eta}}^*}
 & H^q_c(\mathcal{X}_\red,R\Psi_{\!\mathcal{X},c}\Lambda)\ar[d]^-{f^*}\\
 H^q_c\bigl(t(\mathcal{Y})_{\overline{\eta}},\Lambda\bigr)\ar[r]^-{\varepsilon_*}
 & H^q_c(\mathcal{Y}_\red,R\Psi_{\!\mathcal{Y},c}\Lambda)\lefteqn{.}
 }
 \]
\end{prop}

\begin{prf}
 Let $\mathcal{X}'\longrightarrow \mathcal{X}$ be an object of $\Phi_{\mathcal{X}}$.
 By Lemma \ref{lem:pseudo-comp-adm-bu} ii), 
 we may take an object $\mathcal{Y}'\longrightarrow \mathcal{Y}$ of $\Phi_{\mathcal{Y}}$
 such that there exists a (unique) pseudo-compactifiable morphism 
 $f'\colon \mathcal{Y'}\longrightarrow \mathcal{X}'$ which makes the following diagram commutative:
 \[
  \xymatrix{%
 \mathcal{Y}'\ar[r]^-{f'}\ar[d]& \mathcal{X}'\ar[d]\\
 \mathcal{Y}\ar[r]^-{f}& \mathcal{X}\lefteqn{.}
 }
 \]
 Let $\mathcal{U}'$ be the inverse image of $\mathcal{X}'_{\varpi}$ under
 $f'\colon \mathcal{Y}'\longrightarrow \mathcal{X}'$. We may define the homomorphism
 $\varepsilon_{f'}\colon H^q_c(t(\mathcal{U}')_{\overline{\eta}},\Lambda)\longrightarrow H^q_c(\mathcal{Y}_\red,R\Psi_{\!\mathcal{Y},c}\Lambda)$ in the same way as in Definition \ref{defn:e-X}. 
 By Lemma \ref{lem:xi-property} ii), the following diagram is commutative:
 \[
 \xymatrix{%
 H^q_c\bigl(t(\mathcal{X}'_\varpi)_{\overline{\eta}},\Lambda\bigr)\ar[r]^-{\varepsilon_{\mathcal{X}'}}\ar[d]^-{t(f')_{\overline{\eta}}^*}& H^q_c(\mathcal{X}_\red,R\Psi_{\!\mathcal{X},c}\Lambda)\ar[d]^-{f^*}\\
 H^q_c\bigl(t(\mathcal{U}')_{\overline{\eta}},\Lambda\bigr)\ar[r]^-{\varepsilon_{f'}}& H^q_c(\mathcal{Y}_\red,R\Psi_{\!\mathcal{Y},c}\Lambda)\lefteqn{.}
 }
 \]
 Moreover, by Lemma \ref{lem:xi-property} iii), $\varepsilon_{f'}$ is equal to the composite of the push-forward
 $H^q_c(t(\mathcal{U}')_{\overline{\eta}},\Lambda)\longrightarrow H^q_c(t(\mathcal{Y})_{\overline{\eta}},\Lambda)$
 and $\varepsilon_*$. Therefore, the two homomorphism $f^*\circ \varepsilon_*$ and $\varepsilon_*\circ t(f)_{\overline{\eta}}^*$ are equal on the image of $H^q_c(t(\mathcal{X}'_\varpi)_{\overline{\eta}},\Lambda)\longrightarrow H^q_c(t(\mathcal{X})_{\overline{\eta}},\Lambda)$. By Corollary \ref{cor:H_c-limit},
 we have $f^*\circ \varepsilon_*=\varepsilon_*\circ t(f)_{\overline{\eta}}^*$, as desired.
\end{prf}

\subsubsection{Comparison result}
\begin{thm}\label{thm:comp-H_c}
 Assume that a formal scheme $\mathcal{X}$ over $\mathcal{S}$ satisfies the following conditions:
 \begin{itemize}
  \item quasi-compact and separated over $\mathcal{S}$,
  \item pseudo-compactifiable,
  \item and locally algebraizable.
 \end{itemize}
Then, $\varepsilon_*\colon H^q_c(t(\mathcal{X})_{\overline{\eta}},\Lambda)\longrightarrow H^q_c(\mathcal{X}_\red,R\Psi_{\!\mathcal{X},c}\Lambda)$ is an isomorphism.

 In particular, if $\mathcal{X}$ is quasi-compact and algebraizable, 
 then $\varepsilon_*$ can be defined and is an isomorphism.
\end{thm}

First we reduce to the case where $\mathcal{X}$ is algebraizable:

\begin{lem}\label{lem:Cech}
 To prove Theorem \ref{thm:comp-H_c}, we may assume that $\mathcal{X}$ is algebraizable.
\end{lem}

\begin{prf}
 Assume that Theorem \ref{thm:comp-H_c} is true if the formal scheme is algebraizable.
 Let $\mathcal{X}=\bigcup_{i\in I}\mathcal{U}_i$ be an open covering of $\mathcal{X}$,
 where $\mathcal{U}_i$ is algebraizable. Then we have the following morphism of spectral sequences
 (\cf \cite[Remark 5.5.12 iii)]{MR1734903}, Lemma \ref{lem:xi-property} iv)), where we put
 $\mathcal{U}_{i_0}\cap \cdots\cap \mathcal{U}_{i_p}=\mathcal{U}_{i_0,\ldots,i_p}$:
 \[
  \xymatrix{
 E_1^{-p,q}=\bigoplus_{(i_0,\ldots,i_p)\in I^{p+1}}H^q_c\bigl(t(\mathcal{U}_{i_0,\ldots,i_p})_{\overline{\eta}},\Lambda\bigr)\ar@{=>}[r]\ar[d]^-{\varepsilon_*}& H^{-p+q}_c\bigl(t(\mathcal{X})_{\overline{\eta}},\Lambda\bigr)\ar[d]^-{\varepsilon_*}\\
 E_1^{-p,q}=\bigoplus_{(i_0,\ldots,i_p)\in I^{p+1}}H^q_c\bigl((\mathcal{U}_{i_0,\ldots,i_p})_\red,R\Psi_{\mathcal{U}_{i_0,\ldots,i_p},c}\Lambda\bigr)\ar@{=>}[r]& 
 H^{-p+q}_c(\mathcal{X}_\red,R\Psi_{\!\mathcal{X},c}\Lambda)\lefteqn{.}
 }
 \]
 Since $\mathcal{U}_{i_0}\cap \cdots\cap \mathcal{U}_{i_p}=\mathcal{U}_{i_0,\ldots,i_p}$ is algebraizable (\cf \cite[Lemma 7.1.4 (ii)]{MR2309992}), the morphism $\varepsilon_*$ on $E_1$-term is an isomorphism. Therefore
 $\varepsilon_*$ for $\mathcal{X}$ is also an isomorphism.
\end{prf}

Let $X$ be a separated $S$-scheme of finite type, $Y$ a reduced closed subscheme of $X_s$
and $i\colon Y\hooklongrightarrow X_s$ the natural closed immersion.
In the sequel, we consider the case where $\mathcal{X}$ is the completion of $X$ along $Y$.

\begin{defn}
 Let $e\colon (t(\mathcal{X})_{\overline{\eta}})_{\et}\longrightarrow (X_{\overline{\eta}})_\et$
 denotes the morphism of sites constructed similarly as $\varepsilon$ in the diagram 
 in Lemma \ref{lem:natural-morphism}.
 For $L\in D^+(X_{\overline{\eta}},\Lambda)$, let
 $e^*\colon i^*R\psi_X L\longrightarrow R\Psi^\ad_{\!\mathcal{X}} e^*L$ be the morphism constructed in the same way as
 $\varepsilon^*$ in Theorem \ref{thm:comp-Berkovich}.
 If $L\in D^+_c(X_{\overline{\eta}},\Lambda)$, then $e^*$ is an isomorphism by \cite[Theorem 3.1]{MR1395723}.
 We also denote by $e^*$ the morphism induced on the cohomology
 \[
  e^*\colon H^q(Y,i^*R\psi_X L)\longrightarrow H^q(\mathcal{X}_\red,R\Psi^\ad_{\!\mathcal{X}}e^*L)=H^q\bigl(t(\mathcal{X})_{\overline{\eta}},e^*L\bigr).
 \]
\end{defn}

Let $\Phi_X$ denotes the category of blow-ups of $X$ whose centers are contained in $Y$.
Then $\Phi_X$ is cofiltered and the natural functor from $\Phi_X$ to $\Phi_{\mathcal{X}}$ is cofinal.
For an object $X'\longrightarrow X$ of $\Phi_X$, we denote the corresponding object of $\Phi_{\mathcal{X}}$
by $\mathcal{X}'\longrightarrow \mathcal{X}$ and take an open subscheme $X'_\varpi$ of $X'$
such that the $\varpi$-adic completion of $X'_\varpi$ coincides with $\mathcal{X}'_\varpi$. 
We denote the natural morphism of sites 
$(t(\mathcal{X}'_\varpi)_{\overline{\eta}})_{\et}\longrightarrow ((X'_\varpi)_{\overline{\eta}})_\et$
by $e'$.

\begin{defn}
 For $L\in D^+_c(X_{\overline{\eta}},\Lambda)$, we put $L'=L\vert_{(X'_\varpi)_{\overline{\eta}}}$
 and define the homomorphism 
 $e_{X'}\colon H^q_c(t(\mathcal{X}'_\varpi)_{\overline{\eta}},(e^*L)\vert_{t(\mathcal{X}'_\varpi)_{\overline{\eta}}})=H^q_c(t(\mathcal{X}'_\varpi)_{\overline{\eta}},e'^*L')\longrightarrow H^q_c(Y,Ri^!R\psi_X L)$
  as the composite of
 \begin{align*}
 H^q_c\bigl(t(\mathcal{X}'_\varpi)_{\overline{\eta}},e'^*L'\bigr)&\yrightarrow[\cong]{\xi^{-1}_{\mathcal{X}'}}
 H^q_c\bigl((\mathcal{X}'_\varpi)_\red,R\Psi^\ad_{\!\mathcal{X}'_\varpi} e'^*L'\bigr)\yleftarrow[\cong]{e'^*}
 H^q_c\bigl((X'_\varpi)_s,R\psi_{X'_\varpi} L'\bigr)\\
 &\longrightarrow H^q_c(Y,Ri^!R\psi_X L).
 \end{align*}
 Obviously $e_{X'}$ is independent of the choice of $X'_\varpi$ and compatible with the change of $X'$.
 Therefore, by taking inductive limit with respect to $\Phi_X^{\mathrm{op}}$, we obtain the homomorphism
 \[
  e_*\colon H^q_c\bigl(t(\mathcal{X})_{\overline{\eta}},e^*L\bigr)\longrightarrow H^q_c(Y,Ri^!R\psi_X L).
 \]
\end{defn}
It is easy to see that the composite of $e_*\colon H^q_c(t(\mathcal{X})_{\overline{\eta}},\Lambda)\longrightarrow H^q_c(Y,Ri^!R\psi_X\Lambda)$ and the isomorphism
$H^q_c(Y,Ri^!R\psi_X\Lambda)\yrightarrow{\cong}H^q_c(\mathcal{X}_\red,R\Psi_{\!\mathcal{X},c}\Lambda)$ in 
Proposition \ref{prop:comparison-scheme} is equal to $\varepsilon_*$.
Thus it suffices to show that $e_*$ is an isomorphism. The reason why we prefer $e^*$, $e_*$ 
is that they allow a coefficient $L$.

In order to prove the isomorphy of $e_*$, we use the duality theory. 
First we prepare a comparison result for dualizing sheaves.

\begin{prop}\label{prop:comparison-dual}
 Consider the following diagram, where $X_{\overline{\eta}}^\ad=t(\mathcal{S})_{\overline{\eta}}\times_{\Spec \overline{F}}X_{\overline{\eta}}$:
 \[
  \xymatrix{%
 t(\mathcal{X})_{\overline{\eta}}\ar@/^20pt/[rr]^-{e}\ar[r]^-{j}\ar[rd]_-{t(\varphi)_{\overline{\eta}}}&
 X_{\overline{\eta}}^\ad\ar[d]^-{\varphi^\ad}\ar[r]^-{e_1}& X_{\overline{\eta}}\ar[d]^-{\varphi}\\
 & t(\mathcal{S})_{\overline{\eta}}\ar[r]^-{e}& \Spec \overline{F}\lefteqn{.}
 }
 \]
 For $L\in D^+_c(\Spec \overline{F},\Lambda)$, we have natural isomorphisms
 \[
  e_1^*R\varphi^! L\yrightarrow{\cong}R\varphi^{\ad!} e^*L,\qquad 
  e^*R\varphi^! L\yrightarrow{\cong}Rt(\varphi)_{\overline{\eta}}^! e^*L.
 \]
\end{prop}

For this proposition, first we give the following:

\begin{prop}\label{prop:comparison-bikkuri}
 Let $f\colon Y\longrightarrow X$ be a separated morphism between $\overline{F}$-schemes of finite type.
 Let $f^\ad\colon Y^\ad\longrightarrow X^\ad$ be the induced morphism of adic spaces,
 where $(-)^\ad=t(\mathcal{S})_{\overline{\eta}}\times_{\Spec \overline{F}}(-)$. 
 Denote the natural morphisms of sites $(X^\ad)_\et\longrightarrow X_\et$ and 
 $(Y^\ad)_\et\longrightarrow Y_\et$ by $\varepsilon$.
 Then for $L\in D^+(X,\Lambda)$ we have a natural morphism
 $\varepsilon^*Rf^! L\yrightarrow{\cong}Rf^{\ad!} \varepsilon^*L$.
 Moreover, if $L\in D^+_c(X,\Lambda)$, then it is an isomorphism.
\end{prop}

\begin{prf}
 First we will construct a morphism $\varepsilon^*Rf^! L\longrightarrow Rf^{\ad!} \varepsilon^*L$.
 By \cite[Theorem 5.7.5]{MR1734903}, we have the isomorphism of functors 
 $\varepsilon^*\circ Rf_!\yrightarrow{\cong} Rf^\ad_!\circ \varepsilon^*$. This induces a morphism
 $Rf^\ad_!\varepsilon^*Rf^! L\yleftarrow{\cong} \varepsilon^*Rf_!Rf^!L\yrightarrow{\adj}\varepsilon^*L$,
 which corresponds to a morphism
 $\varepsilon^*Rf^! L\longrightarrow Rf^{\ad!}\varepsilon^*L$.

 Next we will prove that this is an isomorphism if $L\in D^+_c(X,\Lambda)$. Since we may work locally on $X$
 and $Y$,
 we may assume that there exists a factorization $Y\hooklongrightarrow\A^d_X\longrightarrow X$ of $f$,
 where the first morphism is a closed immersion. Thus it suffices to consider the case where $f$ is
 the natural morphism $\A^d_X\longrightarrow X$ or the case where $f$ is a closed immersion.

 Consider the first case. By the construction of the trace
 map for adic spaces \cite[proof of Theorem 7.3.4]{MR1734903}, we have the following commutative diagram:
 \[
 \xymatrix{%
 \varepsilon^*Rf_!f^* L(d)[2d]\ar[rr]^-{\varepsilon^*\Tr_f}\ar[d]^-{\cong}&& \varepsilon^*L\ar@{=}[dd]\\
 Rf^\ad_!\varepsilon^*f^* L(d)[2d]\ar@{=}[d]\\
 Rf^\ad_!f^{\ad *}\varepsilon^*L(d)[2d]\ar[rr]^-{\Tr_{f^\ad}}&& \varepsilon^*L\lefteqn{.}
 }
 \]
 This gives the following commutative diagram, whose horizontal arrows are isomorphism by
 \cite[Theorem 7.5.3]{MR1734903}:
 \[
  \xymatrix{%
 \varepsilon^*f^*L(d)[2d]\ar[r]^-{\cong}\ar@{=}[d]& \varepsilon^*Rf^!L\ar[d]\\
 f^{\ad *}\varepsilon^*L(d)[2d]\ar[r]^-{\cong}& Rf^{\ad!}\varepsilon^*L\lefteqn{.}
 }
 \]
 Therefore the right vertical arrow is an isomorphism, as desired.

 Consider the second case. Let $j\colon X\setminus Y\hooklongrightarrow X$ be the natural open immersion.
 Then we have the following morphism of distinguished triangles:
 \[
  \xymatrix{%
 \varepsilon^*Rf^!L\ar[r]\ar[d]& \varepsilon^*f^*L\ar[r]\ar@{=}[d] & \varepsilon^*f^*Rj_*j^*L\ar[r]\ar[d]^-{(*)}&\varepsilon^*Rf^!L[1]\ar[d]\\
 Rf^{\ad!}\varepsilon^*L\ar[r]& f^{\ad*}\varepsilon^*L\ar[r] & f^{\ad*}Rj^\ad_*Rj^{\ad*}\varepsilon^*L\ar[r]&
 Rf^{\ad!}\varepsilon^*L[1]\lefteqn{.}
 }
 \]
 Therefore it suffice to show that $(*)$ is an isomorphism, which is a consequence of
 \cite[Theorem 3.8.1]{MR1734903} since $j^*L$ is constructible.
\end{prf}

\begin{prf}[of Proposition \ref{prop:comparison-dual}]
 By Proposition \ref{prop:comparison-bikkuri}, we have an isomorphism
 $e_1^*R\varphi^! L\yrightarrow{\cong} R\varphi^{\ad!}e^*L$.
 Since $j$ is an open immersion (\cf \cite[Proposition 1.9.6]{MR1734903}, \cite[Lemma 3.13 i)]{MR1620118}),
 by taking $j^*$ of it, we get an isomorphism $e^*R\varphi^! L\yrightarrow{\cong}Rt(\varphi)_{\overline{\eta}}^! e^*L$. 
\end{prf}

Put $K_{t(\mathcal{X})_{\overline{\eta}}}=Rt(\varphi)_{\overline{\eta}}^!\Lambda$,
$K_{X_{\overline{\eta}}^\ad}=R\varphi^{\ad!}\Lambda$ and 
$K_{X_{\overline{\eta}}}=\varphi^!\Lambda$. By Proposition \ref{prop:comparison-dual},
we have $e_1^*K_{X_{\overline{\eta}}}\cong K_{X_{\overline{\eta}}^{\ad}}$ and
$e^*K_{X_{\overline{\eta}}}\cong K_{t(\mathcal{X})_{\overline{\eta}}}$.
Moreover, by the construction of these isomorphisms, the following diagram is commutative:
\[
 \xymatrix{%
 H^0_c\bigl(t(\mathcal{X})_{\overline{\eta}},K_{t(\mathcal{X})_{\overline{\eta}}}\bigr)\ar[rd]_-{\Adj_{t(\mathcal{X})_{\overline{\eta}}}}\ar[r]
 & H^0_c\bigl(X^\ad_{\overline{\eta}},K_{X^\ad_{\overline{\eta}}}\bigr)\ar[d]^-{\Adj_{X^\ad_{\overline{\eta}}}}
 & H^0_c\bigl(X_{\overline{\eta}},K_{X_{\overline{\eta}}}\bigr)\ar[ld]^-{\Adj_{X_{\overline{\eta}}}}
 \ar[l]^-{\cong}_-{e_1^*}\\
 &\Lambda\lefteqn{.}
 }
\]
In the diagram above, $\Adj_\bullet$ denote the natural adjunction homomorphisms.

We have natural cup products
\begin{align*}
 H^p\bigl(t(\mathcal{X})_{\overline{\eta}},K_{t(\mathcal{X})_{\overline{\eta}}}\bigr)\otimes H^q_c\bigl(t(\mathcal{X})_{\overline{\eta}},\Lambda\bigr)&\yrightarrow{\cup} H^{p+q}_c\bigl(t(\mathcal{X})_{\overline{\eta}},K_{t(\mathcal{X})_{\overline{\eta}}}\bigr),\\
% H^p(\mathcal{X}_\red,R\Psi^\ad K_{t(\mathcal{X})_{\overline{\eta}}})\otimes H^q_c\bigl(\mathcal{X}_\red,R\Psi^\ad\Lambda\bigr)&\yrightarrow{\cup} H^{p+q}_c(\mathcal{X}_\red,R\Psi^\ad K_{t(\mathcal{X})_{\overline{\eta}}}),\\
 H^p(Y,i^*R\psi_X K_{X_{\overline{\eta}}})\otimes H^q_c(Y,Ri^!R\psi_X\Lambda)&\yrightarrow{\cup} H^{p+q}_c(Y,Ri^!R\psi_X K_{X_{\overline{\eta}}}).
\end{align*}

\begin{lem}\label{lem:projection-formula}
 We have $e_*(e^*(x)\cup y)=x\cup e_*(y)$ for every 
 $x\in H^p(Y,i^*R\psi_X K_{X_{\overline{\eta}}})$ and 
 $y\in H^q_c(t(\mathcal{X})_{\overline{\eta}},\Lambda)$.
 Note that since $K_{X_{\overline{\eta}}}$ is constructible,
 \[
  e_*\colon H^{p+q}_c(t(\mathcal{X})_{\overline{\eta}},K_{t(\mathcal{X})_{\overline{\eta}}})
 \yleftarrow{\cong} H^{p+q}_c(t(\mathcal{X})_{\overline{\eta}},e^*K_{X_{\overline{\eta}}})
 \yrightarrow{e_*} H^{p+q}_c(Y,Ri^!R\psi_X K_{X_{\overline{\eta}}})
 \]
 can be defined.
\end{lem}

\begin{prf}
 Let $X'\longrightarrow X$ be an object of $\Phi_X$ and $\mathcal{X}'\longrightarrow \mathcal{X}$
 the corresponding object of $\Phi_{\mathcal{X}}$.
 We may assume that $y$ lies in the image of $H^q_c(t(\mathcal{X}'_\varpi)_{\overline{\eta}},\Lambda)\longrightarrow H^q_c(t(\mathcal{X})_{\overline{\eta}},\Lambda)$; thus it suffice to show the analogous property
 for the diagram below:
 \[
  \xymatrix{%
  H^p\bigl((\mathcal{X}'_\varpi)_\red,R\Psi^\ad_{\!\mathcal{X}'_{\varpi}} K_{t(\mathcal{X}'_{\varpi})_{\overline{\eta}}}\bigr)\otimes H^q_c\bigl((\mathcal{X}'_\varpi)_\red,R\Psi^\ad_{\!\mathcal{X}'_\varpi} \Lambda\bigr)\ar@<60pt>[d]^-{\xi_{\mathcal{X}'_{\varpi}}}\ar[r]^-{\cup}\ar@<-40pt>@{=}[d]
 & H^{p+q}_c\bigl((\mathcal{X}'_\varpi)_\red,R\Psi^\ad_{\!\mathcal{X}'_\varpi} K_{t(\mathcal{X}'_{\varpi})_{\overline{\eta}}}\bigr)\ar[d]^-{\xi_{\mathcal{X}'_\varpi}}\\
 H^p\bigl(t(\mathcal{X}'_\varpi)_{\overline{\eta}},K_{t(\mathcal{X}'_{\varpi})_{\overline{\eta}}}\bigr)\otimes H^q_c\bigl(t(\mathcal{X}'_\varpi)_{\overline{\eta}},\Lambda\bigr)\ar[r]^-{\cup}& H^{p+q}_c\bigl(t(\mathcal{X}'_\varpi)_{\overline{\eta}},K_{t(\mathcal{X}'_{\varpi})_{\overline{\eta}}}\bigr)\lefteqn{.}
 }
 \]
 It easily follows from the definition of $\xi_{\mathcal{X}'_\varpi}$.
\end{prf}

\begin{prop}\label{prop:trace-compare}
 Let $\Adj_X\colon H^0_c(Y,Ri^!R\psi_X K_{X_{\overline\eta}})\longrightarrow \Lambda$ be the composite of
 \[
  H^0_c(Y,Ri^!R\psi_X K_{X_{\overline\eta}})\longrightarrow H^0_c(Y,Ri^!K_{X_s})=H^0_c(Y,K_Y)
 \yrightarrow{\adj} \Lambda,
 \]
 where $K_{X_s}$ and $K_Y$ denote the dualizing sheaf of $X_s$ and $K_Y$, respectively.
 Then we have $\Adj_X\circ e_*=\Adj_{t(\mathcal{X})_{\overline{\eta}}}$.
\end{prop}

\begin{prf}
 Let $X'\longrightarrow X$ be an object of $\Phi_X$.
 By the commutative diagram
 \[
  \xymatrix{%
 H^0_c\bigl(t(\mathcal{X}'_\varpi)_{\overline{\eta}},K_{t(\mathcal{X}'_\varpi)_{\overline{\eta}}}\bigr)\ar[r]\ar[d]^-{e'_*}\ar@/^20pt/[rrr]^-{\Adj_{t(\mathcal{X}'_{\varpi})_{\overline{\eta}}}}& 
 H^0_c\bigl(t(\mathcal{X})_{\overline{\eta}},K_{t(\mathcal{X})_{\overline{\eta}}}\bigr)\ar[d]^-{e_*}\ar[rr]_-{\Adj_{t(\mathcal{X})_{\overline{\eta}}}}&& \Lambda\ar@{=}[d]\\
 H^0_c\bigl((X'_\varpi)_s,R\psi_{X'_\varpi} K_{(X'_\varpi)_{\overline{\eta}}}\bigr)\ar[r]\ar@/_20pt/[rrr]_-{\Adj_{X'_\varpi}}&
 H^0_c\bigl(Y,Ri^!R\psi_{X'_\varpi} K_{X_{\overline{\eta}}}\bigr)\ar[rr]^-{\Adj_X}&& \Lambda\lefteqn{,}
 }
 \]
 we may replace $X$ by $X'_\varpi$. Namely, we may assume that $Y=X_s$.

 Let us consider the natural morphism $H^0_c(X_s,R\psi_X K_{X_{\overline{\eta}}})\longrightarrow H^0_c(X_{\overline{\eta}},K_{X_{\overline{\eta}}})$. Then the following diagram is commutative:
 \[
  \xymatrix{%
 H^0_c(\mathcal{X}_\red,R\Psi^\ad_{\!\mathcal{X}} K_{t(\mathcal{X})_{\overline{\eta}}})\ar[d]^-{\xi_{\mathcal{X}}}_-{\cong}&& H^0_c(X_s,R\psi_X K_{X_{\overline{\eta}}})\ar[d]\ar[ll]_-{e^*}^-{\cong}\\
 H^0_c\bigl(t(\mathcal{X})_{\overline{\eta}},K_{t(\mathcal{X})_{\overline{\eta}}}\bigr)\ar[r]\ar[rd]_-{\Adj_{t(\mathcal{X})_{\overline{\eta}}}}& 
 H^0_c(X^\ad_{\overline{\eta}},K_{X^\ad_{\overline{\eta}}})\ar[d]^-{\Adj_{X_{\overline{\eta}}^\ad}}&
 H^0_c(X_{\overline{\eta}},K_{X_{\overline{\eta}}})\ar[ld]^-{\Adj_{X_{\overline{\eta}}}}\ar[l]_-{e_1^*}^-{\cong}\\
 & \Lambda\lefteqn{.}
 }
 \]
 We have already seen the commutativity of the lower triangles.
 The commutativity of the upper rectangle can be checked easily by using adjointness many times; 
 we use a compactification of $X$ over $S$, since the definition of the isomorphism $\xi_{\mathcal{X}}$
 requires a pseudo-compactification of $\mathcal{X}$. 
 Details are left to the reader.

 Since the composite of
 $H^0_c(X_s,R\psi_XK_{X_{\overline{\eta}}})\longrightarrow H^0_c(X_{\overline{\eta}},K_{X_{\overline{\eta}}})\yrightarrow{\Adj_{X_{\overline{\eta}}}}\Lambda$ is $\Adj_X$,
 the proposition follows immediately from the diagram above.
\end{prf}

\begin{prf}[Theorem \ref{thm:comp-H_c}]
 We remain in the setting introduced after Lemma \ref{lem:Cech}.
 By \cite[Theorem 7.1.1]{MR1734903},
 \[
  H^{-q}\bigl(t(\mathcal{X})_{\overline{\eta}},K_{t(\mathcal{X})_{\overline{\eta}}}\bigr)\otimes
 H_c^q\bigl(t(\mathcal{X})_{\overline{\eta}},\Lambda\bigr)\longrightarrow H_c^0\bigl(t(\mathcal{X})_{\overline{\eta}},K_{t(\mathcal{X})_{\overline{\eta}}}\bigr)\yrightarrow{\Adj_{t(\mathcal{X})_{\overline{\eta}}}}\Lambda
 \]
 gives a perfect pairing.
 On the other hand, by \cite[Th\'eor\`eme 4.2]{MR1293970}, 
 \[
  H^{-q}(Y,i^*R\psi_X K_{X_{\overline{\eta}}})\otimes H^q_c(Y,Ri^!R\psi_X\Lambda)
 \longrightarrow H^0_c(Y,Ri^!R\psi_X K_{X_{\overline{\eta}}})\yrightarrow{\Adj_X}\Lambda
 \]
 also gives a perfect pairing. By Lemma \ref{lem:projection-formula} and Proposition \ref{prop:trace-compare},
 $e_*$ is the transpose of $e^*\colon H^{-q}(Y,i^*R\psi_X K_{X_{\overline{\eta}}})\longrightarrow H^{-q}(t(\mathcal{X})_{\overline{\eta}},K_{t(\mathcal{X})_{\overline{\eta}}})$
 under these pairings. Since $e^*$ is an isomorphism, $e_*$ is also an isomorphism.
\end{prf}

Finally we prove the $\ell$-adic version of Theorem \ref{thm:comp-H_c}:

\begin{cor}\label{cor:comp-H_c-l-adic}
 In the setting of Theorem \ref{thm:comp-H_c}, assume that the characteristic of $F$ is zero.
 Then we have natural isomorphisms
 \begin{align*}
  \varepsilon_*&\colon H^q_c\bigl(t(\mathcal{X})_{\overline{\eta}},\Z_\ell\bigr)\longrightarrow H^q_c(\mathcal{X}_\red,R\Psi_{\!\mathcal{X},c}\Z_\ell),\\
  \varepsilon_*&\colon H^q_c\bigl(t(\mathcal{X})_{\overline{\eta}},\Q_\ell\bigr)\longrightarrow H^q_c(\mathcal{X}_\red,R\Psi_{\!\mathcal{X},c}\Q_\ell)\lefteqn{.}
 \end{align*}
\end{cor}

\begin{prf}
 It suffices to prove that 
 \begin{align*}
 H^q_c\bigl(t(\mathcal{X})_{\overline{\eta}},\Z_\ell\bigr)&\cong \varprojlim_n H^q_c\bigl(t(\mathcal{X})_{\overline{\eta}},\Z/\ell^n\Z\bigr),\\
 H^q_c(\mathcal{X}_\red,R\Psi_{\!\mathcal{X},c}\Z_\ell)&\cong \varprojlim_n H^q_c(\mathcal{X}_\red,R\Psi_{\!\mathcal{X},c}\Z/\ell^n\Z).
 \end{align*}
 We may assume that $\mathcal{X}$ is algebraizable.
 Then the former follows from \cite[Theorem 3.3 (ii)]{MR1626021} and \cite[Lemma 3.13 i)]{MR1620118}.
 The latter follows from Proposition \ref{prop:comparison-scheme} and its $\ell$-adic version;
 indeed, if $\mathcal{X}$ is the completion of a separated $S$-scheme $X$ of finite type
 along a closed subscheme $Y$ of $X_s$ and
 denote the closed immersion $Y\hooklongrightarrow X$ by $i$, then we have
 \begin{align*}
 \varprojlim_n H^q_c(\mathcal{X}_\red,R\Psi_{\!\mathcal{X},c}\Z/\ell^n\Z)&\cong 
 \varprojlim_n H^q_c(Y,Ri^!R\psi_X\Z/\ell^n\Z)\cong H^q_c(Y,Ri^!R\psi_X\Z_\ell)\\
 &\cong H^q_c(\mathcal{X}_\red,R\Psi_{\!\mathcal{X},c}\Z_\ell).
 \end{align*}
\end{prf}

\def\cprime{$'$}
\providecommand{\bysame}{\leavevmode\hbox to3em{\hrulefill}\thinspace}
\providecommand{\MR}{\relax\ifhmode\unskip\space\fi MR }
% \MRhref is called by the amsart/book/proc definition of \MR.
\providecommand{\MRhref}[2]{%
  \href{http://www.ams.org/mathscinet-getitem?mr=#1}{#2}
}
\providecommand{\href}[2]{#2}

\end{document}